\documentclass[aap,MSNbibl,seceqn,dvips]{arximspdf}
\usepackage{graphicx}
\doi{10.1214/13-AAP965} 
\volume{24}
\issue{5}
\pubyear{2014}
\firstpage{1850}
\lastpage{1888}

\makeatletter
\newcommand{\widebar}{\overline}
\renewcommand{\bigwedge}{\wedge}
\newcommand{\rrvert}{\vert}
\newcommand{\llvert}{\vert}
\newtheorem{teo}{Theorem}
\newproclaim{ass}{Assumption}
\newproclaim{rem}{Remark}
\makeatother

\begin{document}
\begin{frontmatter}

\title{Limit theorems for the empirical distribution function of scaled
increments of It\^o semimartingales at high frequencies\thanksref{T1}}
\runtitle{\hspace*{-12pt}Empirical c.d.f. of scaled increments of It\^o semimartingales}

\begin{aug}
\author[A]{\fnms{Viktor} \snm{Todorov}\corref{}\ead[label=e1]{v-todorov@northwestern.edu}}
\and
\author[B]{\fnms{George} \snm{Tauchen}\ead[label=e2]{george.tauchen@duke.edu}}
\runauthor{V. Todorov and G. Tauchen}
\affiliation{Northwestern University and Duke University}
\address[A]{Department of Finance\\
Northwestern University\\
Evanston, Illinois 60208-2001\\
USA\\
\printead{e1}} 
\address[B]{Department of Economics\\
Duke University\\
Durham, North Carolina 27708-0097\\
USA\\
\printead{e2}}
\end{aug}
\thankstext{T1}{Supported in part by NSF Grant SES-0957330.}

\received{\smonth{2} \syear{2013}}
\revised{\smonth{8} \syear{2013}}

%
\begin{abstract}
We derive limit theorems for the empirical distribution function of
``devolatilized'' increments of an It\^o semimartingale observed at
high frequencies. These ``devolatilized'' increments are formed by
suitably rescaling and truncating the raw increments to remove the
effects of stochastic volatility and ``large'' jumps. We derive the
limit of the empirical c.d.f. of the adjusted increments for any It\^o
semimartingale whose dominant component at high frequencies has
activity index of $1< \beta\le2$, where $\beta= 2$ corresponds to
diffusion. We further derive an associated CLT in the jump-diffusion
case. We use the developed limit theory to construct a feasible and
pivotal test for the class of It\^o semimartingales with nonvanishing
diffusion coefficient against It\^o semimartingales with no diffusion component.
\end{abstract}

%
\begin{keyword}[class=AMS]
\kwd[Primary ]{62F12}
\kwd{62M05}
\kwd[; secondary ]{60H10}
\kwd{60J75}
\end{keyword}
\begin{keyword}
\kwd{It\^o semimartingale}
\kwd{Kolmogorov--Smirnov test}
\kwd{high-frequency data}
\kwd{stochastic volatility}
\kwd{jumps}
\kwd{stable process}
\end{keyword}
\pdfkeywords{62F12, 62M05, 60H10, 60J75, Ito semimartingale,
 Kolmogorov-Smirnov test, high-frequency data, stochastic volatility, jumps, stable process}
\end{frontmatter}

\section{Introduction}\label{secintro}
The standard jump-diffusion model used for modeling many stochastic
processes is an It\^o semimartingale given by the following
differential equation:
%
%
\begin{equation}
\label{eqXcont} dX_t = \alpha_t\,dt+\sigma_{t}\,dW_t+dY_t,
\end{equation}
where $\alpha_t$ and $\sigma_t$ are processes with c\`{a}dl\`{a}g
paths, $W_t$ is a Brownian motion and $Y_t$ is an It\^o semimartingale
process of pure-jump type (i.e., semimartingale with zero second
characteristic, Definition II.2.6 in \cite{JS}).

At high frequencies, provided $\sigma_t$ does not vanish, the dominant
component of~$X_t$ is its continuous martingale component and at these
frequencies the increments of $X_t$ in (\ref{eqXcont}) behave like
scaled and independent Gaussian random variables. That is, for each
fixed $t$, we have the following convergence:
%
%
\begin{equation}
\label{eqlg} \qquad\frac{1}{\sqrt{h}}(X_{t+sh}-X_t)\stackrel{\mathcal {L}} {\longrightarrow }\sigma_t\times(B_{t+s}-B_{s})
\qquad\mbox{as }h\rightarrow0\mbox{ and }s\in[0,1],
\end{equation}
where $B_t$ is a Brownian motion, and the above convergence is for the
Skorokhod topology; see, for example, Lemma~1 of \cite{TT12}. There are
two distinctive features of the convergence in (\ref{eqlg}). The first
is the scaling factor of the increments on the left-hand side of (\ref
{eqlg}) is the square-root of the length of the high-frequency
interval, a~feature that has been used in developing tests for presence
of diffusion. The second distinctive feature is that the limiting
distribution of the (scaled) increments on the right-hand side of (\ref
{eqlg}) is mixed Gaussian (the mixing given by $\sigma_t$). Both these
features of the local Gaussianity result in (\ref{eqlg}) for models in
(\ref{eqXcont}) have been key in the construction of essentially all
nonparametric estimators of functionals of volatility. Examples include
the jump-robust bipower variation of \cite{BNS04a,BNS06} and the many
other alternative measures of powers of volatility summarized in the
recent book of \cite{JP}. Another important example is the general
approach of \cite{MZ08} (see also~\cite{MSS12}) where estimators of
functions of volatility are formed by utilizing directly (\ref{eqlg})
and working as if volatility is constant over a block of decreasing length.

Despite the generality of the jump-diffusion model in (\ref
{eqXcont}), however, there are several examples of stochastic
processes considered in various applications that are not nested in the
model in (\ref{eqXcont}). Examples include pure-jump It\^o
semimartingales [i.e., the model in (\ref{eqXcont}) with $\sigma_t =
0$ and jumps present], semimartingales contaminated with noise or more
generally nonsemimartingales. In all these cases, both the scaling
constant on the left-hand side of (\ref{eqlg}) as well as the limiting
process on the right-hand side of (\ref{eqlg}) change. Our goal in
this paper, therefore, is to derive a limit theory for a feasible
version of the local Gaussianity result in (\ref{eqlg}) based on
high-frequency record of $X$. An application of the developed limit
theory is a feasible and pivotal test based on Kolmogorov--Smirnov type
distance for the class of It\^o semimartingales with nonvanishing
diffusion component.

The result in (\ref{eqlg}) implies that the high-frequency increments
are approximately Gaussian, but the key obstacle of testing directly
(\ref{eqlg}) is that the (conditional) variance of the increments,
$\sigma_t^2$, is unknown and further is approximately constant only
over a short interval of time. Therefore, on a first step we split the
high-frequency increments into blocks (with length that shrinks
asymptotically to zero as we sample more frequently) and form local
estimators of volatility over the blocks. We then scale the
high-frequency increments within each of the blocks by our local
estimates of the volatility. This makes the scaled high-frequency
increments approximately i.i.d. centered normal random variables with
unit variance. To purge further the effect of ``big'' jumps, we then
discard the increments that exceed a time-varying threshold (that
shrinks to zero asymptotically) with time-variation determined by our
estimator of the local volatility. We derive a (functional) central
limit theorem (CLT) for the convergence of the empirical c.d.f. of the
scaled high-frequency increments, not exceeding the threshold, to the
c.d.f. of a standard normal random variable. The rate of convergence can
be made arbitrary close to~$\sqrt{n}$, by appropriately choosing the
rate of increase of the block size, where $n$ is the number of
high-frequency observations within the time interval. This is achieved
despite the use of the block estimators of volatility, each of which
can estimate the spot volatility $\sigma_t$ at a rate no faster than $n^{1/4}$.

We further derive the limit behavior of the empirical c.d.f. described
above in two possible alternatives to the model (\ref{eqXcont}). The
first is the case where $X_t$ does not contain a diffusive component,
that is, the second term in (\ref{eqXcont}) is absent. Models of
these type have received a lot of attention in various fields; see, for
example, \cite{BNS01,Davis,Mikosch,Kluppelberg} and
\cite{Woerner07}. The second alternative to (\ref{eqXcont}) is the
case in which the It\^o semimartingale is distorted with measurement
error. In each of these two cases, the empirical c.d.f. of the scaled
high-frequency increments below the threshold converges to a c.d.f. of a
distribution different from the standard normal law. This is the stable
distribution in the pure-jump case and the distribution of the noise in
the case of It\^o semimartingale observed with error.

The paper is organized as follows. In Section~\ref{secsetup} we
introduce the formal setup and state the assumptions needed for our
theoretical results. In Section~\ref{sececdf} we construct our
statistic and in Sections~\ref{seccp} and~\ref{secclt} we derive its
limit behavior. In Section~\ref{sectv} we construct the statistic using
alternative local estimator of volatility and derive its limit behavior
in the jump-diffusion case. Section~\ref{secksg} constructs a feasible
test for local Gaussianity using our limit theory, and in Sections~\ref{secmc} and \ref{secempirical} we apply the test on simulated and real
financial data, respectively. The proofs are given in Section~\ref{secapp}.

\section{Setup}\label{secsetup}
We start with the formal setup and assumptions. We will generalize the
setup in (\ref{eqXcont}) to accommodate also the alternative
hypothesis in which $X$ can be of pure-jump type. Thus, the generalized
setup we consider is the following. The process $X$ is defined on a
filtered space $(\Omega,\mathcal{F},(\mathcal{F}_{t})_{t\geq
0},\mathbb
{P})$ and has the following dynamics:
%
%
\begin{equation}
\label{eqX} dX_t = \alpha_t\,dt+\sigma_{t-}\,dS_t+dY_t,
\end{equation}
where $\alpha_t$, $\sigma_t$ and $Y_t$ are processes with c\`{a}dl\`
{a}g paths adapted to the filtration, and $Y_t$ is of pure-jump type.
$S_t$ is a stable process with a characteristic function (see, e.g.,
\cite{SATO}), given by
%
%
\begin{eqnarray}
\label{eqstable} \log \bigl[\mathbb{E}\bigl(e^{iuS_t}\bigr) \bigr] &=& -
t|cu|^{\beta} \bigl(1-i\gamma \operatorname{sign}(u)\Phi \bigr),
\nonumber
\\[-8pt]
\\[-8pt]
\Phi &=& \cases{ \displaystyle \tan(\pi\beta/2), &\quad if $\beta\neq1$,
\vspace*{5pt}\cr
\displaystyle -\frac{2}{\pi}\log|u|, &\quad if $\beta=1$,}
\nonumber
\end{eqnarray}
where $\beta\in(0,2]$ and $\gamma\in[-1,1]$. When $\beta=2$ and $c=1/2$
in (\ref{eqstable}), we recover our original jump-diffusion
specification in (\ref{eqXcont}) in the \hyperref[secintro]{Introduction}. When $\beta<2$,
$X$ is of pure-jump type. $Y_t$ in (\ref{eqX}) will play the role of a
``residual'' jump component at high frequencies (see Assumption~\hyperref[assA2]{A2}
below). We note that $Y_t$ can have dependence with $S_t$ ($\alpha
_t$ and $\sigma_t$), and thus $X_t$ does not ``inherit'' the tail
properties of the stable process $S_t$; for example, $X_t$ can be
driven by a tempered stable process whose tail behavior is very
different from that of the stable process.

Throughout the paper we will be interested in the process $X$ over an
interval of fixed length, and hence without loss of generality we will
fix this interval to be $[0,1]$. We collect our basic assumption on the
components in $X$ next.

{\renewcommand{\theass}{A}
\begin{ass}\label{assA}
$X_t$ satisfies (\ref{eqX}).
\begin{longlist}[A2.]
\item[A1.]\label{assA1} $|\sigma_t|^{-1}$ and $|\sigma
_{t-}|^{-1}$ are strictly positive on $[0,1]$. Further, there is a
sequence of stopping times $T_p$ increasing to infinity and for each
$p$ a bounded process $\sigma_t^{(p)}$ satisfying
$t<T_p \Longrightarrow \sigma_t=\sigma_t^{(p)}$ and a positive constant $K_p$ such that
%
%
\begin{equation}
\label{eqa1} \mathbb{E} \bigl(\bigl|\sigma_t^{(p)}-
\sigma_s^{(p)}\bigr|^2|\mathcal {F}_s
\bigr)\leq K_p |t-s|\qquad\mbox{for every } 0\leq s\leq t\leq1.
\end{equation}

\item[A2.]\label{assA2} There is a sequence of stopping times
$T_p$ increasing to infinity and for each~$p$ a process $Y_t^{(p)}$
satisfying $t<T_p \Longrightarrow Y_t=Y_t^{(p)}$ and a positive
constant $K_p$ such that
%
%
\begin{equation}
\label{eqa2} \mathbb{E} \bigl(\bigl|Y_t^{(p)}-Y_s^{(p)}\bigr|^q|
\mathcal{F}_s \bigr)\leq K_p |t-s|\qquad\mbox{for every } 0\leq s\leq t\leq1
\end{equation}
and for every $q> \beta'$ where $\beta'<\beta$.
\end{longlist}
\end{ass}}%

The assumption in (\ref{eqa1}) can be easily verified for It\^o
semimartingales which is the typical way of modeling $\sigma_t$, but it
is also satisfied for models outside of this class. The condition in
(\ref{eqa2}) can be easily verified for pure-jump It\^o
semimartingales; see, for example, Corollary 2.1.9 of \cite{JP}.

\begin{rem}
Our setup in (\ref{eqX}) (together with
Assumption~\ref{assA}) includes the more parsimonious pure-jump models for $X$
of the form $\int_0^t\sigma_{s-}\,dL_s$ and $L_{T_t}$ where $T_t$ is
absolute continuous time-change process, and $L_t$ is a L\'{e}vy
process with no diffusion component and L\'{e}vy density of the form
$\frac{A_+1_{\{x>0\}}+A_-1_{\{x<0\}}}{|x|^{1+\beta}}+\nu'(x)$ for
$|\nu
'(x)|\leq\frac{K}{|x|^{1+\beta'}}$ when $|x|<x_0$ for some $x_0>0$ (and\vspace*{1pt}
assumptions for $\sigma_t$ and the density of the time change as in Assumption~\hyperref[assA1]{A1}
above). We refer to \cite{TT12} and their supplementary appendix where
this is shown.
\end{rem}

Under Assumption~\ref{assA}, we can extend the local Gaussianity result in (\ref
{eqlg}) to
%
%
\begin{equation}
\label{eqls} \qquad h^{-1/\beta}(X_{t+sh}-X_t) \stackrel{\mathcal{L}} {\longrightarrow }\sigma_{t}\times\bigl(S'_{t+s}-S'_t
\bigr)\qquad\mbox{as }h\rightarrow0\mbox{ and }s\in[0,1]
\end{equation}
for every $t$ and where $S'_t$ is a L\'{e}vy process identically
distributed to $S_t$ and the convergence in (\ref{eqls}) being for the
Skorokhod topology; see, for example, Lemma~1 of \cite{TT12}. That is,
the local behavior of the increments of the process is like that of a
stable process in the more general setting of (\ref{eqX}).

For deriving the CLT for our statistic [in the case of the
jump-diffusion model in (\ref{eqXcont})], we need a stronger
assumption which we state next.

{\renewcommand{\theass}{B}
\begin{ass}\label{assB}
$X_t$ satisfies (\ref{eqX})
with $\beta=2$, that is, $S_t = W_t$.
\begin{longlist}[B3.]
\item[B1.]\label{assB1} The process $Y_t$ is of the form
%
%
\begin{equation}
\label{eqb1} Y_t = \int_0^t\!\int
_{E}\delta^Y(s,x)\mu(ds,dx),
\end{equation}
where $\mu$ is Poisson measure on $\mathbb{R}_+\times E$ with L\'{e}vy
measure $\nu(dx)$ and $\delta^Y(t,x)$ is some predictable function on
$\Omega\times\mathbb{R}_+\times E$.

\item[B2.]\label{assB2} $|\sigma_t|^{-1}$ and $|\sigma
_{t-}|^{-1}$ are strictly positive on $[0,1]$. Further, $\sigma_t$ is
an It\^o semimartingale having the following representation:
%
%
\begin{eqnarray}\label{eqb2}
\sigma_t &=& \sigma_0+\int
_0^t\tilde{\alpha}_u\,du+\int
_0^t\tilde{\sigma}_u\,dW_u+ \int_0^t\tilde{\sigma}'_u\,dW'_u
\nonumber\\[-8pt]\\[-8pt]
&&{} +
\int_0^t\!\int_{E}\delta
^{\sigma}(s,x)\mu(ds,dx),\nonumber
\end{eqnarray}
where $W'_t$ is a Brownian motion independent from $W_t$; $\tilde{\alpha}_t$, $\tilde{\sigma}_t$ and $\tilde{\sigma}'_t$ are
processes with c\`{a}dl\`{a}g paths and $\delta^{\sigma}(t,x)$ is a
predictable function on $\Omega\times\mathbb{R}_+\times E$.

\item[B3.]\label{assB3} $\tilde{\sigma}_t$ and
$\tilde{\sigma}'_t$ are It\^o semimartingales with coefficients with c\`
{a}dl\`
{a}g paths and further jumps being integrals of some predictable
functions, $\delta^{\tilde{\sigma}}$ and $\delta^{\tilde{\sigma}'}$, with respect to the jump measure $\mu$.

\item[B4.]\label{assB4} There is a sequence of stopping times
$T_p$ increasing to infinity and for each~$p$ a deterministic
nonnegative function $\gamma_p(x)$ on $\mathbb{E}$, satisfying $\nu(x\dvtx
\gamma_p(x)\neq0)<\infty$ and such that $| \delta^{Y}(t,x)|\wedge1+|
\delta^{\sigma}(t,x)|\wedge1+| \delta^{\tilde{\sigma}}(t,x)|\wedge
1+| \delta^{\tilde{\sigma}'}(t,x)|\wedge1\leq\gamma_p(x)$ for
$t\leq T_p$.
\end{longlist}
\end{ass}}%

The It\^o semimartingale restriction on $\sigma_t$ (and its
coefficients) is satisfied in most applications. Similarly, we allow
for general time-dependence in the jumps in $X$ which encompasses most
cases in the literature. Assumption~\hyperref[assB4]{B4} is the strongest assumption, and it
requires the jumps to be of finite activity.

\section{\texorpdfstring{Empirical\hspace*{1.5pt} CDF\hspace*{1.5pt} of\hspace*{1.5pt} the\hspace*{1.5pt} ``devolatilized''\hspace*{1.5pt} high-frequency\hspace*{1.5pt} increments}
{Empirical CDF of the ``devolatilized'' high-frequency increments}}\label{sececdf}
Throughout the paper we assume that $X$ is observed on the equidistant
grid $0,\frac{1}{n},\ldots,1$ with $n\rightarrow\infty$. In the
derivation of our statistic we will suppose that $S_t$ is a Brownian
motion and then in the next section we will derive its behavior under
the more general case when $S_t$ is a stable process. The result in
(\ref{eqlg}) suggests that the high-frequency increments $\Delta_i^nX
= X_{i/n}-X_{(i-1)/n}$ are approximately Gaussian with
conditional variance given by the value of the process $\sigma_t^2$ at
the beginning of the increment. Of course, the stochastic volatility
$\sigma_t$ is not known and varies over time. Hence to test for the
local Gaussianity of the high-frequency increments we first need to
estimate locally $\sigma_t$ and then divide the high-frequency
increments by this estimate. To this end, we divide the interval
$[0,1]$ into blocks each of which contains $k_n$ increments, for some
deterministic sequence $k_n\rightarrow\infty$ with $k_n/n\rightarrow
0$. On each of the blocks our local estimator of $\sigma_t^2$ is given
by\looseness=1
%
%
\begin{equation}
\label{eqvhat} \qquad\widehat{V}_j^n = \frac{\pi}{2}
\frac{n}{k_n-1}\sum_{i=(j-1)k_n+2}^{jk_n}\bigl|
\Delta_{i-1}^nX\bigr|\bigl|\Delta_i^nX\bigr|,\qquad j=1,\ldots,\lfloor n/k_n\rfloor.
\end{equation}\looseness=0
$\widehat{V}_j^n$ is the bipower variation proposed by \cite
{BNS04a,BNS06} for measuring the quadratic variation of the diffusion
component of $X$. We note that an alternative measure of $\sigma_t$ can
be constructed using the so-called truncated variation. It turns out,
however, that while the behavior of the two volatility measures in the
case of the jump-diffusion model (\ref{eqXcont}) is the same, it
differs in the case when $S_t$ is stable with $\beta<2$. Using
truncated variation will lead to degenerate limit of our statistic,
unlike the case of using the bipower variation estimator in (\ref
{eqvhat}). For this reason we prefer the latter in our analysis, but
later in Section~\ref{sectv} we also derive in the jump-diffusion case
the behavior of the statistic when truncated variation is used in its
construction.\looseness=1

We use the first $m_n$ increments on each block, with $m_n\leq k_n$, to
test for local Gaussianity. The case $m_n = k_n$ amounts to using all
increments in the block and we will need $m_n<k_n$ for deriving
feasible CLT-s later on. Finally, we remove the high-frequency
increments that contain ``big'' jumps. The total number of increments
used in our statistic is thus given by
%
%
\begin{equation}
\label{eqnjcount} N^n(\alpha,\varpi) = \sum
_{j=1}^{\lfloor n/k_n\rfloor} \sum_{i=(j-1)k_n+1}^{(j-1)k_n+m_n}1
\Bigl(\bigl|\Delta_i^n X\bigr|\leq\alpha\sqrt {
\widehat{V}{}^n_{j}}n^{-\varpi} \Bigr),
\end{equation}
where $\alpha>0$ an $\varpi\in(0,1/2)$. We note that here we use a
time-varying threshold in our truncation to account for the
time-varying $\sigma_t$.

The scaling of every high-frequency increment will be done after
adjusting $\widehat{V}_j^n$ to exclude the contribution of that
increment in its formation
%
%
\begin{equation}
\label{eqvhati} \qquad \widehat{V}_j^n(i) = \cases{
\displaystyle \frac{k_n-1}{k_n-3}\widehat{V}_j^n -
\frac{\pi}{2}\frac{n}{k_n-3}\bigl|\Delta_{i}^nX\bigr|\bigl|
\Delta_{i+1}^nX\bigr|,
\vspace*{3pt}\cr
\qquad\mbox{for }i=(j-1)k_n+1,
\vspace*{6pt}
\cr
\displaystyle \frac{k_n-1}{k_n-3}\widehat{V}_j^n
- \frac{\pi}{2}\frac{n}{k_n-3} \bigl(\bigl|\Delta_{i-1}^nX\bigr|\bigl|
\Delta_{i}^nX\bigr|+\bigl|\Delta_{i}^nX\bigr|\bigl|
\Delta_{i+1}^nX\bigr| \bigr),
\vspace*{3pt}\cr
\qquad\mbox{for }i=(j-1)k_n+2,\ldots,jk_n-1,
\vspace*{6pt}
\cr
\displaystyle \frac{k_n-1}{k_n-3}
\widehat{V}_j^n - \frac{\pi}{2}\frac{n}{k_n-3}\bigl|
\Delta_{i-1}^nX\bigr|\bigl|\Delta_{i}^nX\bigr|,
\vspace*{3pt}\cr
\qquad\mbox{for }i=jk_n.}
\end{equation}
With this, we define
%
%
\begin{eqnarray}\label{eqfhat}
\widehat{F}_n(\tau) &=&  \frac{1}{N^n(\alpha,\varpi)}
\nonumber\\[-8pt]\\[-8pt]
&&{}\times \sum_{j=1}^{\lfloor n/k_n\rfloor} \sum_{i=(j-1)k_n+1}^{(j-1)k_n+m_n}1
\biggl\{\frac{\sqrt {n}\Delta_i^nX}{\sqrt{\widehat{V}{}^n_{j}(i)}}\leq\tau \biggr\} 1_{ \{
|\Delta_i^n X|\leq\alpha\sqrt{\widehat{V}{}^n_{j}}n^{-\varpi}\}},\nonumber
\end{eqnarray}
which is simply the empirical c.d.f. of the ``devolatilized'' increments
that do not contain ``big'' jumps. In the jump-diffusion case of (\ref
{eqXcont}), $\widehat{F}_n(\tau)$ should be approximately the c.d.f.
of a
standard normal random variable.

We note that all the results that follow for $\widehat{F}_n(\tau)$ will
continue to hold if we do not truncate for the jumps in the
construction of $\widehat{F}_n(\tau)$.
The intuition for this is easiest to form in the case when $X$ is a L\'
{e}vy process without\vspace*{2pt} drift from the following $\mathbb{E}\llvert 1_{\{
\sqrt{n}\Delta_i^nX\leq\tau\}}-1_{\{\sqrt{n}\sigma\Delta
_i^nW\leq\tau\}
}\rrvert  = O(n^{\beta'/2-1+\iota})$ for $\beta'$ the constant of
Assumption~\hyperref[assA2]{A2} and $\iota>0$ arbitrary small. Our rational for looking
at the truncated increments only is that the order of magnitude of the
above difference; that is, the error due to the presence of jumps in
$X$ can be slightly reduced by using truncation.

The construction of our statistic resembles the practice of
standardizing increments of the process of fixed length by a measure
for volatility constructed from high-frequency data within the interval
(after correcting for jumps and leverage effect); see, for example,
\cite{ABD}. The main difference is that here the length of the
increments that are standardized is shrinking and further the
volatility estimator is local, that is, over a shrinking time interval.
Both these differences are crucial for deriving our feasible limit
theory for~$\widehat{F}_n(\tau)$.

\section{Convergence in probability of \texorpdfstring{$\widehat{F}_n(\tau)$}{Fn(tau)}}\label{seccp}
We next derive the limit behavior of $\widehat{F}_n(\tau)$ both under
the null of model (\ref{eqXcont}) as well as under a set of
alternatives. We start with the case when $X_t$ is given by (\ref{eqX}).
%
%
\begin{teo}\label{teoecdf}
Suppose Assumption~\ref{assA} holds, and assume the block size grows at the rate
%
%
\begin{equation}
\label{eqecdf0} k_n \sim n^{q}\qquad\mbox{for some }q \in(0,1)
\end{equation}
and $m_n\rightarrow\infty$ as $n\rightarrow\infty$. Then if $\beta
\in
(1,2]$, we have
%
%
\begin{equation}
\label{eqecdf1} \widehat{F}_n(\tau)  \stackrel{\mathbb{P}} {\longrightarrow } F_{\beta
}(\tau)\qquad\mbox{as }n\rightarrow\infty,
\end{equation}
where the above convergence is uniform in $\tau$ over compact subsets
of $\mathbb{R}$, $F_{\beta}(\tau)$ is the c.d.f. of $\sqrt{\frac
{2}{\pi
}}\frac{S_1}{\mathbb{E}|S_1|}$ ($S_1$ is the\vspace*{1pt} value of the $\beta
$-stable process $S_t$ at time $1$) and $F_2(\tau)$ equals the c.d.f. of a
standard normal variable $\Phi(\tau)$.
\end{teo}
Since $\widehat{F}_n(\tau)$ and $F_{\beta}(\tau)$ are c\`{a}dl\`{a}g
and nondecreasing, the above result holds also uniformly on $\mathbb{R}$.

\begin{rem}
The limit result in (\ref{eqecdf1})
shows that when $S_t$ is stable with $\beta<2$, $\widehat{F}_n(\tau)$
estimates the c.d.f. of a $\beta$-stable random variable. We note that
when $\beta<2$, the correct scaling factor for the high-frequency\vspace*{1pt}
increments is $n^{1/\beta}$. However, in this case we need also to
scale $\widehat{V}_j^n$ by $n^{1/\beta-1/2}$ in order\vspace*{-1pt} for the latter to
converge to a nondegenerate limit (that is proportional to $\sigma
_t^2$). Hence the ratio
%
%
\begin{equation}
\frac{\sqrt{n}\Delta_i^nX}{\sqrt{\widehat{V}{}^n_{j}(i)}} = \frac
{n^{1/\beta}\Delta_i^nX}{\sqrt{n^{2/\beta-1}\widehat{V}{}^n_{j}(i)}}
\end{equation}
is appropriately scaled even in the case when $\beta<2$ and importantly
without knowing a priori the value of $\beta$. We further note that the
limiting c.d.f., $F_{\beta}(\tau)$, is of a random variable that has the
same scale regardless of the value of $\beta$. That is, in all cases of
$\beta$, $F_{\beta}(\tau)$ corresponds to the c.d.f. of a random variable
$Z$ with $\mathbb{E}|Z| = \sqrt{\frac{2}{\pi}}$. Therefore,\vspace*{1pt} the
difference between $\beta<2$ and the null $\beta=2$ will be in the
relative probability assigned to ``big'' versus ``small'' values of
$\tau$.
\end{rem}

We note further that in Theorem~\ref{teoecdf} we restrict $\beta>1$.
The reason is that for $\beta\leq1$, the limit behavior of $\widehat{F}_n(\tau)$ is determined by the drift term in $X$ (when present) and
not $S_t$. To allow for $\beta\leq1$ and still have a limit result of
the type in~(\ref{eqecdf1}), we need to use $\Delta_i^nX-\Delta
_{i-1}^nX$ in the construction of $\widehat{F}_n(\tau)$ which
essentially eliminates the drift term.

We next derive the limiting behavior of $\widehat{F}_n(\tau)$ in the
situation when the It\^o semimartingale $X$ is ``contaminated'' by
noise, which is of particular relevance in financial applications.
%
%
\begin{teo}\label{teonoise}
Suppose Assumption~\ref{assA} holds and $k_n\propto n^q$ for some $q\in(0,1)$
and $m_n\rightarrow\infty$ as $n\rightarrow\infty$. Let $\widehat
{F}_n(\tau)$ be given by (\ref{eqfhat}) with $\Delta_i^nX$ replaced
with $\Delta_i^nX^*$ for $X^*_{i/n} = X_{i/n}+\varepsilon
_{i/n}$ and where $ \{\varepsilon_{i/n} \}
_{i=1,\ldots,n}$ are i.i.d. random variables defined on a product
extension of the original probability space and independent from
$\mathcal{F}$. Further, suppose $\mathbb{E}|\varepsilon_{i/n}|^{1+\iota}<\infty$ for some $\iota>0$. Finally, \mbox{assume} that the
c.d.f. of $\frac{1}{\mu} (\varepsilon_{i/n}-\varepsilon_{(i-1)/n} )$, $F_{\varepsilon}(\tau)$, is continuous where we denote
$\mu= \sqrt{\frac{\pi}{2}}\sqrt{\mathbb{E} (|\varepsilon
_{i/n}-\varepsilon_{(i-1)/n}||\varepsilon_{(i-1)/n}-\varepsilon
_{(i-2)/n}| )}$. Then
%
%
\begin{equation}
\label{eqnoise1} \widehat{F}_n(\tau) \stackrel{\mathbb{P}} {\longrightarrow } F_{\varepsilon
}(\tau)\qquad\mbox{as }n\rightarrow\infty,
\end{equation}
where the above convergence is uniform in $\tau$ over compact subsets
of $\mathbb{R}$.
\end{teo}

\begin{rem}
When $X$ is observed with noise, the
noise becomes the leading component at high frequencies. Hence, our
statistic recovers the c.d.f. of the (appropriately scaled) noise
component. Similar to the pure-jump alternative of $S_t$ with $\beta
<2$, here $\sqrt{n}$ is not the right scaling for the increments
$\Delta
_i^nX^*$, but this is offset in the ratio in $\widehat{F}_n(\tau)$ by a
scaling factor for the local variance estimator $\widehat{V}_j^n$ that
makes it nondegenerate. Unlike the pure-jump alternative, in the
presence of noise the correct scaling of the numerator and the
denominator in the ratio in $\widehat{F}_n(\tau)$ is given by
%
%
\begin{equation}
\frac{\sqrt{n}\Delta_i^nX^*}{\sqrt{\widehat{V}{}^n_{j}(i)}} = \frac
{\Delta
_i^nX^*}{\sqrt{n^{-1}\widehat{V}{}^n_{j}(i)}},
\end{equation}
that is, we need to scale down $\widehat{V}{}^n_{j}(i)$ to ensure it
converges to nondegenerate limit.
\end{rem}

The limit result in (\ref{eqnoise1}) provides an important insight
into the noise by studying its distribution. We stress the fact that
the presence of $\widehat{V}{}^n_{j}$ in the truncation is very important
for the limit result in (\ref{eqnoise1}). This is because it ensures
that the threshold is ``sufficiently'' big so that it does not matter
in the asymptotic limit. If, on the~other hand, the threshold did not
contain $\widehat{V}{}^n_{j}$ (i.e., $\widehat{V}{}^n_{j}$ was replaced by
$1$ in the threshold), then in this case the limit will be determined
by the behavior of the density of the noise around zero.

We finally note that when $\varepsilon_{i/n}$ is normally
distributed, a case that has received a lot of attention in the
literature, the limiting c.d.f. $F_{\varepsilon}(\tau)$ is that of a centered
normal but with variance that is below $1$. Therefore, in this case
$F_{\varepsilon}(\tau)$ will be below the c.d.f. of a standard normal
variable, $\Phi(\tau)$, when $\tau<0$ and the same relationship will
apply to $1-F_{\varepsilon}(\tau)$ and $1-\Phi(\tau)$ when
$\tau>0$.\vadjust{\goodbreak}

On a more general level, the above results show that the empirical c.d.f.
estimator $\widehat{F}_n(\tau)$ can shed light on the potential sources
of violation of the local Gaussianity of high-frequency data. It
similarly can provide insights on the performance of various estimators
that depend on this hypothesis.

\section{CLT of \texorpdfstring{$\widehat{F}_n(\tau)$}{Fn(tau)} under local Gaussianity}\label{secclt}
%
%
\begin{teo}\label{teosv}
Let $X_t$ satisfy (\ref{eqX}) with $S_t$ being a Brownian motion and
assume that Assumption~\ref{assB} holds. Further, let the block size grow at the rate
%
%
\begin{equation}
\label{eqsv1} \frac{m_n}{k_n}\rightarrow0,\qquad k_n \sim n^{q}\mbox{ for some }q\in(0,1/2)\mbox{ when }n\rightarrow
\infty,
\end{equation}
such that $\frac{nm_n}{k_n^3}\rightarrow\lambda\geq0$. We then have
locally uniformly in subsets of $\mathbb{R}$
%
%
\begin{eqnarray}
\widehat{F}_n(\tau)-\Phi(\tau)
&=&  \widehat{Z}_1^n(\tau)+\widehat {Z}_2^n(
\tau)\nonumber
\\
\label{eqsv2} &&{}
+\frac{1}{k_n} \frac{\tau^2\Phi''(\tau)-\tau\Phi
'(\tau)}{8} \biggl( \biggl(\frac{\pi}{2}\biggr)^2+\pi-3 \biggr)
\\
&&{} +o_p \biggl( \frac{1}{k_n} \biggr)\nonumber
\end{eqnarray}
and
\begin{equation}\label{eqsv3}
\bigl(\matrix{\sqrt{\lfloor n/k_n\rfloor m_n}
\widehat{Z}_1^n(\tau) & \sqrt {\lfloor n/k_n
\rfloor k_n}\widehat{Z}_2^n(\tau) }\bigr)
\stackrel {\mathcal{L}} {\longrightarrow}  \bigl(\matrix{Z_1(
\tau) & Z_2(\tau) }\bigr),
\end{equation}
where $\Phi(\tau)$ is the c.d.f. of a standard normal variable
and $Z_1(\tau)$ and $Z_2(\tau)$ are two independent Gaussian processes
with covariance functions
%
%
\begin{eqnarray}
\label{eqsv4} \operatorname{Cov} \bigl(Z_1(\tau_1),
Z_1(\tau_2) \bigr) &=& \Phi (\tau_1\wedge
\tau _2)-\Phi(\tau_1)\Phi(\tau_2),\nonumber
\\
\operatorname{Cov} \bigl(Z_2(\tau_1), Z_2(
\tau_2) \bigr) &=& \biggl[\frac{\tau_1\Phi
'(\tau_1)}{2} \frac{\tau_2\Phi'(\tau_2)}{2} \biggr]
\biggl( \biggl(\frac{\pi}{2} \biggr)^2+\pi-3 \biggr),
\\
\eqntext{\tau _1,\tau_2\in \mathbb{R}.}
\end{eqnarray}
\end{teo}
Due to the ``big'' jumps, we derive the CLT only on compact sets of
$\tau$ since the error in the estimation of the c.d.f. for $\tau
\rightarrow
\pm\infty$ is affected by the truncation.

We make several observations regarding the limiting result in (\ref
{eqsv2})--(\ref{eqsv4}). The first term of $\widehat{F}_n(\tau
)-\Phi
(\tau)$ in (\ref{eqsv2}), $\widehat{Z}_1^n(\tau)$, converges to
$Z_1(\tau)$ which is the standard Brownian bridge appearing in the
Donsker theorem for empirical processes; see, for example, \cite{Vaart}. The second and third terms on the right-hand side of (\ref
{eqsv2}) are due to the estimation error in recovering the local
variance, that is, the presence of $\widehat{V}_j^n$ in $\widehat
{F}{}^n(\tau)$ instead of the true (unobserved) $\sigma_t^2$. $\widehat
{Z}_2^n(\tau)$ converges to a centered Gaussian process, independent
from $Z_1(\tau)$, while the third term on the right-hand side of (\ref
{eqsv2}) is an asymptotic bias. Importantly, the asymptotic bias as
well as the variance of $ (Z_1(\tau) \ \  Z_2(\tau) )$ are all
constants that depend only on $\tau$ and not the stochastic volatility
$\sigma_t$. Therefore, feasible inference based on (\ref{eqsv2}) is
straightforward.

We note that by picking the rate of growth of $m_n$ and $k_n$
arbitrarily close to~$\sqrt{n}$, we can make the rate of convergence of
$\widehat{F}_n(\tau)$ arbitrary close to $\sqrt{n}$. We should further
point out that this is unlike the rate of estimating the spot $\sigma
_t^2$ by $\widehat{V}_j^n$ (with the same choice of $k_n$) which is at
most $n^{1/4}$. The reason for the better rate of convergence of our
estimator is in the integration of the error due to the estimation~$\widehat{V}_j^n$.

The order of magnitude of the three components on the right-hand side
of (\ref{eqsv2}) are different with the second term always dominated
by the other two. Its presence should provide a better finite-sample
performance of a test based on (\ref{eqsv2}).

Finally, we point out that a feasible CLT for $\widehat{F}_n(\tau)$ is
available with ``only'' arbitrarily close to $\sqrt{n}$ rate of
convergence and not exactly $\sqrt{n}$. This is due to the presence of
the drift term in $X$. The latter leads to asymptotic bias which is of
order $1/\sqrt{n}$ and removing it via de-biasing is in general
impossible as we cannot estimate the latter from high-frequency record
of $X$.

\section{Empirical CDF of ``devolatilized'' high-frequency increments with an alternative volatility estimator}\label{sectv}
As mentioned in Section~\ref{sececdf}, an alternative estimator of the
volatility is the truncated variation of \cite{Ma2} defined as
%
%
\begin{eqnarray}
\label{eqtvhat} \widehat{C}_j^n = \frac{n}{k_n}\sum
_{i=(j-1)k_n+1}^{jk_n}\bigl|\Delta _i^nX\bigr|^21
\bigl(\bigl|\Delta_i^n X\bigr|\leq\alpha n^{-\varpi}
\bigr),\nonumber\\[-12pt]\\[-12pt]
\eqntext{j=1,\ldots,\lfloor n/k_n\rfloor,}
\end{eqnarray}
where $\alpha>0$ and $\varpi\in(0,1/2)$ and the corresponding one
excluding the contribution of the $i$th increment, for
$i=(j-1)k_n+1,\ldots,jk_n$, is
%
%
\begin{equation}
\label{eqtvhati} \widehat{C}_j^n(i) = \frac{k_n}{k_n-1}
\widehat{C}_j^n-\frac
{n}{k_n-1}\bigl|\Delta_i^nX\bigr|^21
\bigl(\bigl|\Delta_i^n X\bigr|\leq\alpha n^{-\varpi
} \bigr).
\end{equation}
We define the corresponding empirical c.d.f. of the ``devolatilized'' (and
truncated) high-frequency increments as
%
%
\begin{equation}
\label{eqfhattv} \widehat{F}_{n}'(\tau) =
\frac{1}{N^{\prime n}(\alpha,\varpi)}\sum_{j=1}^{\lfloor n/k_n\rfloor} \sum
_{i=(j-1)k_n+1}^{(j-1)k_n+m_n}1 \biggl\{ \frac{\sqrt{n}\Delta_i^nX}{\sqrt{\widehat{C}{}^n_{j}(i)}}\leq
\tau \biggr\} 1_{ \{|\Delta_i^n X|\leq\alpha n^{-\varpi} \}},\hspace*{-25pt}
\end{equation}
where for $\alpha>0$ and $\varpi\in(0,1/2)$
%
%
\begin{equation}
\label{eqnjcounttv} N^{\prime n}(\alpha,\varpi) = \sum
_{j=1}^{\lfloor n/k_n\rfloor} \sum_{i=(j-1)k_n+1}^{(j-1)k_n+m_n}1
\bigl(\bigl|\Delta_i^n X\bigr|\leq\alpha n^{-\varpi} \bigr).
\end{equation}
In the next theorem we derive a CLT for $\widehat{F}_{n}'(\tau)$ when
$X$ is a jump-diffusion.\vadjust{\goodbreak}

%
\begin{teo}\label{teotv}
Let $X_t$ satisfy (\ref{eqX}) with $S_t$ being a Brownian motion and
assume that Assumption~\ref{assB} holds. Let $k_n$ and $m_n$ satisfy (\ref
{eqsv1}). We then have locally uniformly in subsets of $\mathbb{R}$
%
%
\begin{equation}
\label{eqtv1} \qquad \widehat{F}'_n(\tau)-\Phi(\tau) =
\widehat{Z}_1^n(\tau)+\widehat {Z}_2^n(
\tau)+\frac{1}{k_n} \frac{\tau^2\Phi''(\tau)-\tau\Phi
'(\tau
)}{4}+o_p \biggl(
\frac{1}{k_n} \biggr)
\end{equation}
and
\begin{equation}\label{eqtv2}
\bigl(\matrix{\sqrt{\lfloor n/k_n\rfloor m_n}
\widehat{Z}_1^n(\tau) & \sqrt {\lfloor n/k_n
\rfloor k_n}\widehat{Z}_2^n(\tau) }\bigr)
\stackrel {\mathcal{L}} {\longrightarrow} \bigl(\matrix{Z_1(
\tau) & Z_2(\tau) }\bigr),
\end{equation}
where $\Phi(\tau)$ is the c.d.f. of a standard normal variable
and $Z_1(\tau)$ and $Z_2(\tau)$ are two independent Gaussian processes
with covariance functions
%
%
\begin{eqnarray}
\label{eqtv3} \operatorname{Cov} \bigl(Z_1(\tau_1),
Z_1(\tau_2) \bigr) &=& \Phi (\tau_1\wedge
\tau _2)-\Phi(\tau_1)\Phi(\tau_2),
\nonumber\\[-8pt]\\[-8pt]
\operatorname{Cov} \bigl(Z_2(\tau_1), Z_2(
\tau_2) \bigr) &=& \bigl[\tau_1\Phi'(
\tau_1) \tau_2\Phi'(\tau_2)
\bigr], \qquad \tau_1,\tau_2\in\mathbb{R}.\nonumber
\end{eqnarray}
Further, in the case when $\alpha_t$, $\sigma_t$ and $\delta^{Y}(t,x)$
do not depend on $t$, the above result continues to hold even when
Assumption~\hyperref[assB4]{\textup{B4}} is replaced with the weaker condition $\int_{E}[|\delta
^{Y}(x)|^{\beta'}\wedge1]\nu(dx)<\infty$ for some $0\leq\beta'<1$,
provided for $\iota>0$ arbitrary small, we have
%
%
\begin{eqnarray}\label{eqtv4}
k_n \bigl(n^{1-(4-\beta')\varpi}\vee n^{-((1-\beta'/2)/(1+\beta
'))+\iota}\vee
n^{-(2/3)(2-\beta')\varpi+\iota} \bigr)&\rightarrow &0,
\nonumber\\[-8pt]\\[-8pt]
\frac{k_n^{3}}{m_n}n^{-(4-\beta')\varpi} &\rightarrow&0.\nonumber
\end{eqnarray}
\end{teo}
The CLT for $\widehat{F}'_n(\tau)$ is similar to that for $\widehat
{F}_n(\tau)$ with the only difference being that the asymptotic bias
[the third term on the right-side of (\ref{eqtv1})] and the limiting
Gaussian process $Z_2$ are of smaller magnitude and with smaller
variance, respectively. This is not surprising as the truncated
variation is known to be a more efficient estimator of volatility than
the bipower variation.

The last part of the theorem shows that in the case when $\alpha_t$,
$\sigma_t$ and $\delta^{Y}(t,x)$ do not depend on $t$, the CLT result
continues to hold in presence of jumps of infinite activity (but finite
variation) provided the growth condition (\ref{eqtv4}) holds. This
condition can be simplified when one uses a value for $\varpi$
arbitrarily close to $1/2$ (as is common) and $m_n$ close to $k_n$.

\section{Test for local Gaussianity of high-frequency data}\label{secksg}
We proceed with a feasible test for a jump-diffusion model of the type
given in (\ref{eqXcont}) using the developed limit theory above. We
will use $\widehat{F}_n(\tau)$ for this. The critical region of our
proposed test is given by
%
%
\begin{equation}
\label{eqcr} C_n = \Bigl\{\sup_{\tau\in\mathcal{A}}
\sqrt{N^n(\alpha,\varpi )}\bigl|\widehat{F}_n(\tau)-\Phi(
\tau)\bigr|> q_n(\alpha,\mathcal{A}) \Bigr\},
\end{equation}
where we recall that $\Phi(\tau)$ denotes the c.d.f. of a standard normal
random variable, $\alpha\in(0,1)$, $\mathcal{A}\in\mathbb{R}$ is a
finite union of compact sets with positive Lebesgue measure and
$q_n(\alpha,\mathcal{A})$ is the $(1-\alpha)$-quantile of
%
%
\begin{eqnarray}
&& \sup_{\tau\in\mathcal{A}} \biggl\llvert Z_1(\tau) + \sqrt{
\frac
{m_n}{k_n}}Z_2(\tau)
\nonumber\\[-8pt]\\[-8pt]
&&\hspace*{18pt}{} + \sqrt{\frac{m_n}{k_n}}
\frac{\sqrt{n}}{k_n} \frac
{\tau^2\Phi''(\tau)-\tau\Phi'(\tau)}{8} \biggl( \biggl(\frac
{\pi
}{2}
\biggr)^2+\pi-3 \biggr)\biggr\rrvert\nonumber
\end{eqnarray}
with $Z_1(\tau)$ and $Z_2(\tau)$ being the Gaussian processes defined
in Theorem~\ref{teosv}. We can easily evaluate $q_n(\alpha,\mathcal
{A})$ via simulation.

We note that in (\ref{eqcr}) we use $N^n(\alpha,\varpi)$ as a
normalizing constant. This is justified because we have $\frac
{N^n(\alpha,\varpi)}{\lfloor n/k_n\rfloor m_n}\stackrel{\mathbb
{P}}{\longrightarrow}1$, both in\vspace*{1pt} the jump-diffusion case as well as in
the two alternative scenarios considered in Section~\ref{seccp}. The
choice of $k_n$ and $m_n$ in general should be dictated by how much
volatility of volatility in $X$ we have. We illustrate this in the next section.

The test in (\ref{eqcr}) resembles a Kolmogorov--Smirnov type test for
equality of continuous one-dimensional distributions. There are two
differences between our test and the original Kolmogorov--Smirnov test.
First, in our test we scale the high-frequency increments by a
nonparametric local estimator of the volatility, and this has an
asymptotic effect on the test statistic, as evident from Theorem~\ref{teosv}. The second difference is in the region $\mathcal{A}$ over
which the difference $\widehat{F}_n(\tau)-\Phi(\tau)$ is evaluated. For
reasons we already discussed, that are particular to our problem here,
we need to exclude arbitrary values of $\tau$ that are high in magnitude.

Now, in terms of the size and power of the test, under Assumptions~\ref{assA}
and~\ref{assB}, using Theorems~\ref{teoecdf}~and~\ref{teosv}, we have
%
%
\begin{eqnarray}
\lim_n\mathbb{P} (C_n ) &=& \alpha\qquad\mbox{if }\beta =2\quad\mbox{and}
\nonumber\\[-8pt]\\[-8pt]
\liminf_n\mathbb{P} (C_n ) &=& 1\qquad\mbox{if }\beta\in(1,2),\nonumber
\end{eqnarray}
where we make also use of the fact that the stable and standard normal
variables have different c.d.f.'s on compact subsets of $\mathbb{R}$ with
positive Lebesgue measure. By Theorem~\ref{teonoise}, the above power
result applies also to the case when we observe $X_{i/n}+\varepsilon_{i/n}$, provided of course the limiting c.d.f. of
the noise in (\ref{eqnoise1}) differs from that of the standard
normal on the set $\mathcal{A}$.

We note that existing tests for presence of diffusive component in $X$
are based only on the scaling factor of the high-frequency increments
on the left-hand side of~(\ref{eqls}). However, the limiting result in
(\ref{eqls}) implies much more. Mainly, the distribution of the
``devolatilized'' increments should be stable (and, in particular,
normal in the jump-diffusion case). Our test in (\ref{eqcr}), unlike
earlier work, incorporates this distribution implication of (\ref
{eqls}) as well.

We finally point out that using Theorem~\ref{teosv}, one should be able
to derive alternative tests for the presence of diffusive component in
$X$, by adopting other measures of discrepancy between distributions
like the Cram\'er--von Mises test.

\section{Monte Carlo}\label{secmc}
We now evaluate the performance of our test on simulated data. We
consider the following two models. The first is
%
%
\begin{eqnarray}
\label{eqmc1} dX_t &=& \sqrt{V_t}\,dW_t+\int
_{\mathbb{R}}x\mu(ds,dx),
\nonumber\\[-8pt]\\[-8pt]
dV_t &=& 0.03(1.0-V_t)\,dt+0.1
\sqrt{V_t}\,dB_t,\nonumber
\end{eqnarray}
where $(W_t,B_t)$ is a vector of Brownian motions with $\operatorname
{Corr}(W_t,B_t)
= -0.5$ and $\mu$ is a homogenous Poisson measure with compensator
$\nu
(dt,dx) = dt\otimes\frac{0.25e^{-|x|/0.4472}}{0.4472}\,dx$ which\vspace*{1pt}
corresponds to double exponential jump process with intensity of $0.5$
(i.e., a jump every second day on average). This model is calibrated to
financial data by setting the means of continuous and jump variation
similar to those found in earlier empirical work. Similarly, we allow
for dependence between $X_t$ and $V_t$, that is, leverage effect. The
second model is given by
%
%
\begin{equation}
\label{eqmc2} X_t = S_{T_t}\qquad\mbox{with } T_t
= \int_0^tV_s\,ds,
\end{equation}
where $S_t$ is a symmetric tempered stable martingale with L\'{e}vy
measure $\frac{0.1089e^{-|x|}}{|x|^{1+1.8}}$, and $V_t$ is the
square-root diffusion given in (\ref{eqmc1}). The process in (\ref
{eqmc2}) is a time-changed tempered stable process. The parameters of
$S_t$ are chosen such that it behaves locally like $1.8$-stable process
and it has variance at time $1$ equal to $1$ [as the model in (\ref
{eqmc1})]. For this process the local Gaussianity does not hold and
hence the behavior of the test on data from the model in (\ref
{eqmc2}) will allow us to investigate~the power of the test. We also
consider another alternative to the jump-diffusion, mainly the case
when the process in (\ref{eqmc2}) is contaminated with i.i.d.
Gaussian noise. The variance of the noise is set to $0.01$ consistent
with empirical evidence in \cite{HL}.

We turn next to the implementation of the test. We apply the test to
one year's worth of simulated data which consists of $252$ days (our
unit of time is one trading\vadjust{\goodbreak} day). We consider two sampling frequencies:
$n=100$ and $n=200$ which correspond to sampling every $5$ and $2$
minutes, respectively, in a typical trading day. We experiment with
1--4 blocks per day. In each block we use $75\%$ or $70\%$ of the
increments in the formation of the test, that is, we set $\lfloor
m_n/k_n \rfloor= 0.75$ for $n=100$ and $\lfloor m_n/k_n \rfloor=
0.70$ for $n=200$. We found very little sensitivity of the test with
respect to the choice of the ratio $m_n/k_n$. For the truncation of the
increments, as typical in the literature, we set $\alpha=3.0$ and
$\varpi=0.49$. Finally, the set $\mathcal{A}$ over which the difference
$\widehat{F}_n(\tau)-F(\tau)$ in our test is evaluated is set to
%
%
\begin{equation}
\label{eqA} \mathcal{A} = \bigl[Q(0.01)\dvtx Q(0.40)\bigr]\cup \bigl[Q(0.60)\dvtx Q(0.99)\bigr],
\end{equation}
where $Q(\alpha)$ is the $\alpha$-quantile of standard normal.

%
\begin{table}[b]
\tabcolsep=14pt
\caption{Monte Carlo results for jump-diffusion model (\protect\ref{eqmc1})}\label{tbmcsize}
\begin{tabular*}{\tablewidth}{@{\extracolsep{\fill}}@{}lccccc@{}}
\hline
& \multicolumn{5}{c@{}}{\textbf{Rejection rate}} \\[-6pt]
& \multicolumn{5}{c@{}}{\hrulefill} \\
\textbf{Nominal size} & \multicolumn{3}{c}{\textbf{Kolmogorov--Smirnov test}} & \multicolumn{2}{c@{}}{\textbf{Power variation based test}}\\
\hline
\multicolumn{6}{@{}c@{}}{Sampling frequency $n=100$}\\
  & $k_n=33$ & $k_n=50$ & $k_n=100$ & $           p=1.0$ & $p=1.5$\\
$\alpha =   1\%$   & $0.0$ & $0.8$ & \phantom{0}$5.6$ & $                1.5$ & $     0.7$ \\
$\alpha =   5\%$   & $0.4$ & $4.3$ & $16.8$ & $                7.7$ & $     5.4$ \\[3pt]
\multicolumn{6}{@{}c@{}}{Sampling frequency $n=200$}\\
  & $k_n=50$ & $k_n=67$ & $k_n=200$ & $           p=1.0$ & $p=1.5$\\
$\alpha =   1\%$   & $0.4$ & $0.9$ & $10.3$ & $                1.4$ & $     1.3$ \\
$\alpha =   5\%$  & $1.2$ & $3.2$ & $32.8$ & $                8.1$ & $     6.7$ \\
\hline
\end{tabular*}
\tabnotetext[]{}{\textit{Note}: For the cases with $n=100$ we set $\lfloor
m_n/k_n\rfloor=0.75$ and for the cases with $n=200$ we set $\lfloor
m_n/k_n \rfloor=0.70$. The power variation test is a one-sided test
based on Theorem~2 in \protect\cite{SJ10} with $k=2$ and cutoff $u_n =
7\hat{\sigma}\Delta_n^{0.49}$ with $\hat{\sigma}$ being an estimate of
volatility over the day using bipower variation.}
\end{table}

The results of the Monte Carlo are reported in Tables~\ref{tbmcsize}--\ref{tbmcpower2}. For the smaller sample size, $n=100$, and
with no blocking at all ($k_n=n$) to account for volatility movements
over the day, there are size distortions most noticeable at the
conventional 5 percent level. With two blocks ($\lfloor n/k_n\rfloor
=2$), size is appropriate, while it is seen to have excellent power in
Tables~\ref{tbmcpower} and \ref{tbmcpower2}. But with three blocks on
$n=100$, there are size distortions because the noisy estimates of
local volatility distort the test. Considering the larger sample size
($n=200$), now with three blocks the test's size is
approximately correct while power is excellent. For larger values of
$k_n$ relative to $n$ ($\lfloor n/k_n\rfloor=1$) the time variation in
volatility over the day coupled with the relatively high precision of
estimating a biased version of local volatility, leads to departures
from Gaussianity of the (small) scaled increments and hence the over-rejections.

%
\begin{table}[t]
\tabcolsep=14pt
\caption{Monte Carlo results for pure-jump model (\protect\ref{eqmc2})}\label{tbmcpower}
\begin{tabular*}{\tablewidth}{@{\extracolsep{\fill}}@{}lccccc@{}}
\hline
& \multicolumn{5}{c@{}}{\textbf{Rejection rate}} \\[-6pt]
& \multicolumn{5}{c@{}}{\hrulefill} \\
\textbf{Nominal size} & \multicolumn{3}{c}{\textbf{Kolmogorov--Smirnov test}} & \multicolumn{2}{c@{}}{\textbf{Power variation based test}}\\
\hline
\multicolumn{6}{@{}c}{Sampling frequency $n=100$}\\
  & $k_n=33$ & $k_n=50$ & $k_n=100$ & $           p=1.0$ & $p=1.5$\\
$\alpha =   1\%$   & \phantom{0}$45.9$ & \phantom{0}$95.2$ & \phantom{0}$99.9$ & $               71.1$ & $    14.2$\\
$\alpha =   5\%$   & \phantom{0}$76.6$ & \phantom{0}$99.6$ & $100.0$ & $               91.3$ & $    38.4$\\[3pt]
\multicolumn{6}{@{}c@{}}{Sampling frequency $n=200$}\\
& $k_n=50$ & $k_n=67$ & $k_n=200$\\
$\alpha= 1\%$ & $100.0$ & $100.0$ & $100.0$ & $     97.7$ &$ 32.1$\\
$\alpha= 5\%$ & $100.0$ & $100.0$ & $100.0$ & $     99.4$ &$ 63.3$\\
\hline
\end{tabular*}
\tabnotetext[]{}{\textit{Note}: Notation as in Table~\ref{tbmcsize}.}\vspace*{-3pt}
\end{table}
%
\begin{table}[b]
\tabcolsep=14pt
\caption{Monte Carlo results for pure-jump model (\protect\ref{eqmc2}) plus noise}\label{tbmcpower2}
\begin{tabular*}{\tablewidth}{@{\extracolsep{\fill}}@{}lccccc@{}}
\hline
& \multicolumn{5}{c@{}}{\textbf{Rejection rate}} \\[-6pt]
& \multicolumn{5}{c@{}}{\hrulefill} \\
\textbf{Nominal size} & \multicolumn{3}{c}{\textbf{Kolmogorov--Smirnov test}} & \multicolumn{2}{c@{}}{\textbf{Power variation based test}}\\
\hline
\multicolumn{6}{@{}c@{}}{Sampling frequency $n=100$}\\
  & $k_n=33$ & $k_n=50$ & $k_n=100$ & $           p=1.0$ & $p=1.5$\\
$\alpha =   1\%$   & \phantom{0}$39.0$ & \phantom{0}$91.4$ & \phantom{0}$99.8$ & $               21.7$ & $     1.6$\\
$\alpha =   5\%$   & \phantom{0}$70.9$ & \phantom{0}$99.6$ & $100.0$ & $               46.2$ & $     5.9$\\
\\[3pt]
\multicolumn{6}{@{}c@{}}{Sampling frequency $n=200$}\\
& $k_n=50$ & $k_n=67$ & $k_n=200$\\
$\alpha= 1\%$ & $100.0$ & $100.0$ & $100.0$ & \phantom{0}$     6.5$ &$  0.0$\\
$\alpha= 5\%$ & $100.0$ & $100.0$ & $100.0$ & $     17.3$ &$  0.0$\\
\hline
\end{tabular*}
\tabnotetext[]{}{\textit{Note}: Notation as in Table~\ref{tbmcsize}.}\vspace*{-3pt}
\end{table}

In Tables~\ref{tbmcsize}--\ref{tbmcpower2} we also report the
performance on the simulated data of a test for presence of Brownian
motion in high-frequency data based on (truncated) power variations
computed on two different frequencies, proposed in \cite{SJ10}; see
also \cite{TT09}. This test, unlike the test proposed here, does not
exploit the distributional implication of the local Gaussianity result
in (\ref{eqlg}). We can see from Table~\ref{tbmcsize} that the test
based on the power variations has reasonable behavior under the null of
presence of a diffusion component in $X$. Table~\ref{tbmcpower} further
shows that for the optimal choice of the power ($p=1$), the test has
slightly lower power against the considered pure-jump alternative in
(\ref{eqmc2}) than the Kolmogorov--Smirnov test (when block size is
chosen optimally).

When the pure-jump model is contaminated with noise, the scaling of the
power variations is similar (for the considered frequencies) to that of
a jump-diffusion model observed without noise. Hence, Table~\ref{tbmcpower2} reveals relatively low power of the test based on the
power variations against the alternative of pure-jump process
contaminated with noise. By contrast, the Kolmogorov--Smirnov test
shows almost no change in performance compared with the alternative
when the pure-jump process is observed without noise (Table~\ref{tbmcpower}). The reason is that the Kolmogorov--Smirnov test
incorporates also the distributional implications of (\ref{eqlg}) and,
under the pure-jump plus noise scenario, the scaled high-frequency
increments have a distribution which is very different from standard normal.

%
\begin{figure}

\includegraphics{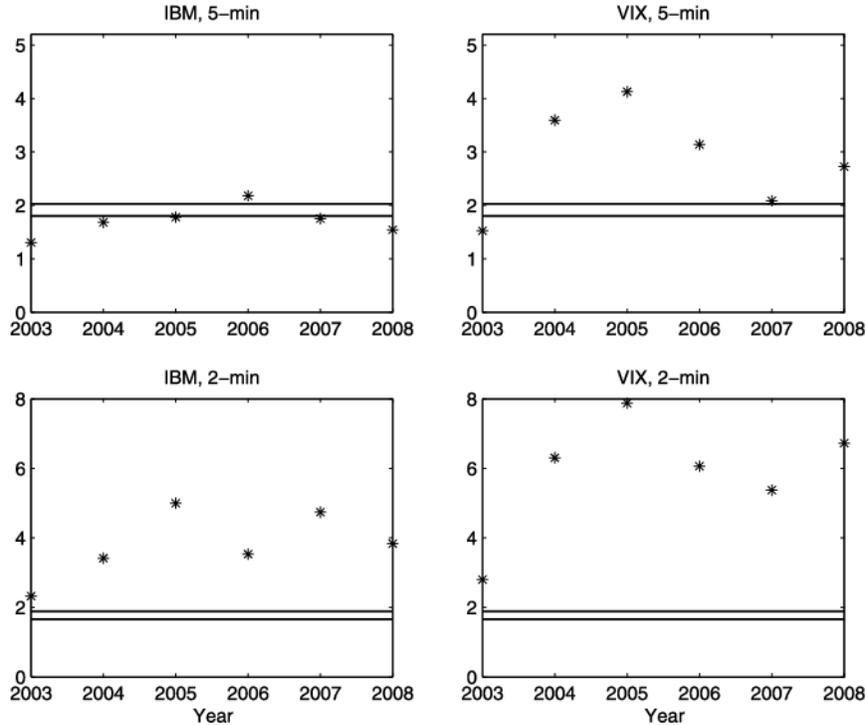}

\caption{Kolmogorov--Smirnov tests for local Gaussianity. The
$*$ corresponds to the value of the test $\sup_{\tau\in A}\sqrt {N^n(\alpha,\varpi)}|\widehat{F}_n(\tau)-F(\tau)|$, and the solid lines
are the critical values $q_n(\alpha,\mathcal{A})$ for $\alpha= 5\%$
and $\alpha= 1\%$.}\label{figtestresults}
\end{figure}

\section{Empirical illustration}\label{secempirical}
We now apply our test to two different financial assets, the IBM stock
price and the VIX volatility index. The analyzed period is
2003--2008, and like in the Monte Carlo we consider two and five
minute sampling frequencies. The test is performed for each of the
years in the sample. We set $\mathcal{A}$ as in (\ref{eqA}) and
$\lfloor n/k_n\rfloor=3$ for the five-minute sampling frequency and
$\lfloor n/k_n\rfloor=4$ for the two-minute frequency. As in the Monte
Carlo, the ratio $\lfloor m_n/k_n\rfloor$ is set to $0.75$ and $0.70$
for the five-minute and two-minute, respectively, sampling frequencies.
Finally, to account for the well-known diurnal pattern in volatility we
standardize the raw high-frequency returns by a time-of-day scale
factor exactly as in \cite{TT10}.

The results from the test are shown on Figure~\ref{figtestresults}. We
can see from the figure that the local Gaussianity hypothesis works
relatively well for the $5$-minute IBM returns. At $2$-minute sampling
frequency for the IBM stock price, however, our test rejects the local
Gaussianity hypothesis at conventional significance levels.
Nevertheless, the values of the test are not very far from the critical
ones. The explanation of the different outcomes of the test on the two
sampling frequencies is to be found in the presence of microstructure
noise. The latter becomes more prominent at the higher frequency.
Turning to the VIX index data, we see a markedly different outcome. For
this data set, the local Gaussianity hypothesis is strongly rejected at
both frequencies. The explanation for this is that the underlying model
is of pure-jump type, that is, the model (\ref{eqX}) with $\beta<2$.

\section{Proofs}\label{secapp}
We start with introducing some notation that we will make use of in the proofs:
%
%
\begin{eqnarray*}
\label{app1} A_t &=& \int_0^t \alpha_s\,ds,\qquad B_t = \int_0^t \sigma_s\,dS_s,\qquad \bar{\sigma}_t =\sigma_t-\sum_{s\leq t}\Delta \sigma_s,
\\
\label{app2} \dot{V}^n_j &=& \frac{n}{k_n-1}\frac{\pi}{2}\sum_{i=(j-1)k_n+2}^{jk_n}\bigl|\Delta_{i-1}^nA+\Delta_{i-1}^nB\bigr|\bigl|\Delta _i^nA+\Delta_{i}^nB\bigr|,
\\
\label{app3} \widetilde{V}{}^n_j &=& \frac{n}{k_n-1}
\frac{\pi}{2}\sum_{i=(j-1)k_n+2}^{jk_n}\bigl|
\Delta_{i-1}^nB\bigr|\bigl|\Delta_i^nB\bigr|,
\\
\widebar{V}_j^n &=& \sigma_{((j-1)k_n)/n}^2
\frac{n}{k_n-1}\frac{\pi}{2}\sum_{i=(j-1)k_n+2}^{jk_n}\bigl|
\Delta_{i-1}^nS\bigr|\bigl|\Delta_i^nS\bigr|
\end{eqnarray*}
and we define $\dot{V}^n_j(i)$, $\widetilde{V}{}^n_j(i)$ and $\widebar{V}{}^n_j(i)$ from the above as in (\ref{eqvhati}). We also denote
%
%
\begin{equation}
\label{app5} \widetilde{F}_n(\tau) = \frac{N^n(\alpha,\varpi)}{\lfloor
n/k_n\rfloor
m_n}
\widehat{F}_n(\tau).
\end{equation}
Finally, in the proofs we will denote with $K$ a positive constant that
might change from line to line but importantly does not depend on $n$
and $\tau$. We will also use the shorthand notation $\mathbb
{E}_i^n
(\cdot ) = \mathbb{E} (\cdot|\mathcal{F}_{(i-1)/n} )$.

\subsection{Localization}

We will prove Theorems~\ref{teoecdf}--\ref{teotv} under the following
stron\-ger versions of Assumption~\ref{assA} and~\ref{assB}:
\begin{longlist}[SA.]
\item[SA.]\label{assSA} We have Assumption~\ref{assA} with $\alpha_t$,
$\sigma_t$ and $\sigma_t^{-1}$ being all uniformly bounded on $[0,1]$.
Further, (\ref{eqa1}) and (\ref{eqa2}) hold for $\sigma_t$ and $Y_t$,
respectively.

\item[SB.]\label{assSB} We have Assumption~\ref{assB} with all processes
$\alpha_t$, $\tilde{\alpha}_t$, $\sigma_t$, $\sigma_t^{-1}$,
$\tilde{\sigma}_t$, $\tilde{\sigma}'_t$ and the
coefficients of
the It\^o semimartingale representations of $\tilde{\sigma}_t$ and
$\tilde{\sigma}'_t$ being uniformly bounded on $[0,1]$. Further
$(|\delta^Y(t,x)| + |\delta^{\sigma}(t,x)| + |\delta^{\tilde{\sigma
}}(t,x)|+ |\delta^{\tilde{\sigma}'}(t,x)|)\leq\gamma(x)$ for some
nonnegative valued function $\gamma(x)$ on $E$ satisfying $\int_E\nu
(x\dvtx  \gamma(x) = 0)\,dx<\infty$ and $\gamma(x)\leq K$ for some constant
$K$.
\end{longlist}

Extending the proofs to the weaker Assumptions~\ref{assA} and \ref{assB} follows by
standard localization techniques exactly as Lemma~4.4.9 of \cite{JP}.

\subsection{Proof of Theorem~\texorpdfstring{\protect\ref{teoecdf}}{1}}
Without loss of
generality, we will assume that $\tau<0$, the case $\tau\geq0$ being
dealt with analogously [by working with $1-\widehat{F}_n(\tau)$ instead].
We first analyze the behavior of $\widehat{V}_j^n$. We denote with
$\eta
_n$ a deterministic sequence that depends only on $n$ and vanishes as
$n\rightarrow\infty$.

Using the triangular inequality, the Chebyshev inequality, successive
conditioning, as well as the H\"{o}lder inequality and Assumption~\hyperref[assSA]{SA},
we get for $j=1,\ldots,\lfloor n/k_n\rfloor$
%
%
\[
\label{appecdf1} \mathbb{P} \bigl(n^{2/\beta-1}\bigl|\widehat{V}_j^n-
\dot{V}_j^n\bigr|\geq \eta_n \bigr)\leq K
\frac{n^{1/\beta-(1/\beta')\wedge1+\iota}}{\eta
_n}\qquad \forall\iota>0.
\]
Similarly, using the triangular inequality, Chebyshev's inequality as
well as the H\"{o}lder inequality, we get for $j=1,\ldots,\lfloor
n/k_n\rfloor$
%
%
\[
\label{appecdf2} \mathbb{P} \bigl(n^{2/\beta-1}\bigl|\dot{V}_j^n-
\widetilde{V}_j^n\bigr|\geq \eta_n \bigr)\leq K
\frac{n^{1/\beta- 1+\iota}}{\eta_n}\qquad \forall\iota >0.
\]
Next, using the triangular inequality, the Chebyshev inequality, the H\"
{o}lder inequality, the Burkholder--Davis--Gundy inequality as well as
Assumption~\hyperref[assSA]{SA}, we get for $j=1,\ldots,\lfloor n/k_n\rfloor$
%
%
\[
\label{appecdf3} \mathbb{P} \bigl(n^{2/\beta-1}\bigl|\widetilde{V}_j^n
- \widebar{V}_j^n\bigr|\geq \eta_n \bigr)\leq
K\frac{k_n^{1/2-\iota}}{n^{1/2-\iota}\eta
_n}\qquad \forall\iota>0.
\]
Finally, using the self-similarity of the stable process and the
Burkholder--Davis--Gundy inequality (for discrete martingales), we get
for $j=1,\ldots,\lfloor n/k_n\rfloor$
%
%
\begin{eqnarray}\label{appecdf4}
\mathbb{P} \biggl(\biggl|n^{2/\beta-1}\widebar{V}_j^n
- \frac{\pi
}{2}\sigma ^2_{((j-1)k_n)/n}\bigl(\mathbb{E}|S_1|\bigr)^2\biggr|
\geq\eta_n \biggr)\leq K \frac{1}{k_n^{\beta-1-\iota}\eta_n^{\beta-\iota}}\nonumber
\\
\eqntext{\forall\iota \in (0,1- \beta).}
\end{eqnarray}
Combining these results, we get altogether for $\forall\iota\in
(0,1-\beta)$
%
%
\begin{eqnarray}
\label{appecdf5} &&\mathbb{P} \biggl(\biggl|n^{2/\beta-1}\widehat{V}_j^n
- \frac{\pi
}{2}\sigma ^2_{((j-1)k_n)/n}\bigl(\mathbb{E}|S_1|\bigr)^2\biggr|
\geq\eta_n \biggr)
\nonumber\\[-8pt]\\[-8pt]
&&\qquad \leq K \biggl(\frac{n^{1/\beta-(1/\beta')\wedge1+\iota
}}{\eta
_n}\vee
\frac{k_n^{1/2-\iota}}{n^{1/2-\iota}\eta_n} \vee \frac
{1}{k_n^{\beta-1-\iota}\eta_n^{\beta-\iota}} \biggr).\nonumber
\end{eqnarray}
Using the same proofs we can show that the result above continues to
hold when $\widehat{V}_j^n$ is replaced with $\widehat{V}_j^n(i)$.

Next, for $i=(j-1)k_n+1,\ldots,(j-1)k_n+m_n$ and $j=1,\ldots,\lfloor
n/k_n\rfloor$, we denote
\[
\cases{ \displaystyle\xi_{i,j}^n(1) =
n^{1/\beta} \biggl( \Delta_i^nA +
\Delta_i^nY +\int_{(i-1)\Delta_n}^{i\Delta_n}(
\sigma_{u-}-\sigma_{((j-1)k_n)/n})\,dS_u \biggr),
\vspace*{5pt}
\cr
\displaystyle\xi_{i,j}^n(2) =
n^{1/\beta}\sigma_{((j-1)k_n)/n}\Delta_i^nS
1_{ \{|\Delta_i^n X|> \alpha\sqrt{\widehat{V}{}^n_{j}}n^{-\varpi } \}}.}
\]
With this notation, using similar inequalities as before, we get
%
%
\begin{equation}
\label{appecdf7} \mathbb{P} \bigl(\bigl|\xi_{i,j}^n(1)\bigr|\geq
\eta_n \bigr)\leq K \biggl(\frac
{n^{1/\beta-(1/\beta')\wedge1+\iota}}{\eta_n}\vee
\frac
{k_n^{\beta
/2+\iota/2}}{n^{\beta/2+\iota/2}\eta_n^{\beta+\iota}} \biggr).
\end{equation}
Next, using the result in (\ref{appecdf5}) above as well as the H\"
{o}lder inequality, we get
%
%
\begin{eqnarray}\label{appecdf8}
&& \mathbb{P} \bigl(\bigl|\xi_{i,j}^n(2)\bigr|\geq \eta_n \bigr)\nonumber
\\
&&\qquad \leq K \frac
{n^{-(1/2-\varpi)\beta+\iota}\vee n^{1/\beta-(1/\beta')\wedge
1+\iota
}\vee(k_n/n)^{1/2-\iota}\vee k_n^{1+\iota-\beta}}{\eta_n^{\iota}}
\\
\eqntext{\forall\iota>0.}
\end{eqnarray}
We next denote the set (note that by Assumption~\hyperref[assSA]{SA}, $\sigma_t$ is
strictly above zero on the time interval $[0,1]$)
%
%
\begin{eqnarray}
\label{appecdf9} \mathcal{A}_{i,j}^n &=& \biggl\{\omega\dvtx
\frac{|\xi_{i,j}^n(1)|+|\xi
_{i,j}^n(2)|}{\sqrt{\pi/2}\mathbb{E}|S_1|}>\eta_n
\nonumber\\[-8pt]\\[-8pt]
&&\hspace*{5pt}{} \cup \biggl\llvert \frac{n^{1/\beta-1/2}\sqrt{\widehat{V}_j^n(i)}}{\sqrt{\pi/2}
\sigma_{((j-1)k_n)/n}\mathbb{E}|S_1|} -1
\biggr\rrvert >\eta_n \biggr\}\nonumber
\end{eqnarray}
for $i=(j-1)k_n+1,\ldots,(j-1)k_n+m_n$ and $j=1,\ldots,\lfloor
n/k_n\rfloor$.

We now can set [recall (\ref{eqecdf0})]
%
%
\begin{equation}
\label{appecdf10} \qquad\eta_n = n^{-x},\qquad 0<x< \biggl[ \biggl(
\frac{1}{\beta'}\wedge1-\frac
{1}{\beta} \biggr)\bigwedge
\frac{1-q}{2} \bigwedge\frac{q(\beta
-1)}{\beta} \biggr]
\end{equation}
and this choice is possible because of the restriction on the rate of
increase of the block size $k_n$ relative to $n$ given in (\ref
{eqecdf0}). With this choice of $\eta_n$, the results in~(\ref{appecdf5}), (\ref{appecdf7}) and (\ref{appecdf8}) imply
%
%
\begin{equation}
\label{appecdf11} \frac{1}{\lfloor n/k_n\rfloor m_n}\sum_{j=1}^{\lfloor n/k_n\rfloor}
\sum_{i=(j-1)k_n+1}^{(j-1)k_n+m_n}\mathbb{P} \bigl(\mathcal
{A}_{i,j}^n \bigr) = o(1).
\end{equation}
Therefore, for any compact subset $\mathcal{A}$ of $(-\infty,0)$,
%
%
\begin{equation}
\label{appecdf12} \sup_{\tau\in\mathcal{A}}\bigl|\widetilde{F}_n(\tau) -
\widehat {G}_n(\tau)\bigr| = o_p(1),
\end{equation}
where we denote
%
%
\begin{eqnarray*}\label{appecdf13}
\widehat{G}_n(\tau) &=& \frac{1}{\lfloor n/k_n\rfloor m_n}
\\
&&{}\times \sum_{j=1}^{\lfloor n/k_n\rfloor} \sum_{i=(j-1)k_n+1}^{(j-1)k_n+m_n}1
\biggl\{ \frac{\sqrt{n}\Delta_i^nX}{\sqrt{\widehat{V}{}^n_{j}(i)}}1_{ \{
|\Delta
_i^n X|\leq\alpha\sqrt{\widehat{V}{}^n_{j}}n^{-\varpi} \}}\leq \tau \biggr\}1_{\{(\mathcal{A}_{i,j}^n)^c\}}.
\end{eqnarray*}
Taking into account the definition of the set $\mathcal{A}_{i,j}^n$,
we get
\[
\cases{ \displaystyle \widehat{G}_n(\tau)\geq
\frac{1}{\lfloor n/k_n\rfloor m_n}\sum_{j=1}^{\lfloor n/k_n\rfloor} \sum
_{i=(j-1)k_n+1}^{(j-1)k_n+m_n}1 \biggl\{\frac{n^{1/\beta}\Delta_i^nS}{\sqrt{\pi/2}\mathbb{E}|S_1|}\leq
\tau(1-\eta_n)-\eta_n \biggr\},
\vspace*{6pt}\cr
\displaystyle\widehat{G}_n(\tau)\leq\frac{1}{\lfloor n/k_n\rfloor m_n}\sum
_{j=1}^{\lfloor n/k_n\rfloor} \sum_{i=(j-1)k_n+1}^{(j-1)k_n+m_n}1
\biggl\{\frac{n^{1/\beta}\Delta_i^nS}{\sqrt{\pi/2}\mathbb{E}|S_1|}\leq\tau(1+\eta_n)+\eta_n
\biggr\}.}
\]

Using the Glivenko--Cantelli theorem (see, e.g., Theorem 19.1 of \cite{Vaart}), we have
%
%
\begin{eqnarray*}
\label{appecdf15} &&\sup_{\tau} \biggl|\frac{1}{\lfloor n/k_n\rfloor m_n}\sum
_{j=1}^{\lfloor
n/k_n\rfloor} \sum_{i=(j-1)k_n+1}^{(j-1)k_n+m_n}1
\biggl\{\frac
{n^{1/\beta
}\Delta_i^nS}{\sqrt{\pi/2}\mathbb{E}|S_1|}\leq\tau (1-\eta _n)-\eta_n
\biggr\}
\\
&&\qquad\hspace*{174pt}{}  - F_{\beta}\bigl(\tau(1-\eta _n)-\eta
_n\bigr) \Biggr|\stackrel{\mathbb{P}} {\longrightarrow}0,
\nonumber
\\
\label{appecdf15a} &&\sup_{\tau} \Biggl|\frac{1}{\lfloor n/k_n\rfloor m_n}\sum
_{j=1}^{\lfloor
n/k_n\rfloor} \sum_{i=(j-1)k_n+1}^{(j-1)k_n+m_n}1
\biggl\{\frac
{n^{1/\beta
}\Delta_i^nS}{\sqrt{\pi/2}\mathbb{E}|S_1|}\leq\tau (1+\eta _n)+\eta_n
\biggr\}
\\
&&\qquad\hspace*{177.5pt}{}- F_{\beta}(\tau(1+\eta _n)+\eta
_n \Biggr|\stackrel{\mathbb{P}} {\longrightarrow}0\nonumber
\end{eqnarray*}
and further using the smoothness of c.d.f. of the stable distribution we have
%
%
\begin{eqnarray*}
\label{appecdf16} \sup_{\tau}\bigl\llvert F_{\beta}\bigl(
\tau(1-\eta_n)-\eta_n\bigr)-F_{\beta
}(\tau)\bigr
\rrvert &\rightarrow& 0,
\\
\sup_{\tau}\bigl\llvert
F_{\beta}\bigl(\tau(1+\eta _n)+\eta _n
\bigr)-F_{\beta}(\tau)\bigr\rrvert &\rightarrow &0.
\end{eqnarray*}
These two results together imply
%
%
\[
\label{appecdf17} \sup_{\tau}\bigl|\widehat{G}_n(
\tau)-F_{\beta}(\tau)\bigr|\stackrel {\mathbb {P}} {\longrightarrow}0\nonumber
\]
and from here, using (\ref{appecdf12}), we have $\sup_{\tau\in
\mathcal
{A}}|\widetilde{F}_n(\tau)-F_{\beta}(\tau)| = o_p(1)$ for any compact
subset $\mathcal{A}$ of $(-\infty,0)$. Hence, to prove (\ref
{eqecdf1}), we need only to show
%
%
\begin{equation}
\label{appecdf19} \frac{N^n(\alpha,\varpi)}{\lfloor n/k_n\rfloor m_n}\stackrel {\mathbb {P}} {\longrightarrow}1\qquad\mbox{as }n\rightarrow\infty.
\end{equation}
We have
%
%
\begin{eqnarray*}
\label{appecdf20}
&& \mathbb{P} \Bigl(\bigl|\Delta_i^n X\bigr|> \alpha
\sqrt{\widehat {V}{}^n_{j}}n^{-\varpi} \Bigr)
\\
&&\qquad \leq \mathbb{P} \biggl(\biggl\llvert \frac
{n^{1/\beta
-1/2}\sqrt{\widehat{V}_j^n}}{\sqrt{\pi/2}\sigma_{((j-1)k_n)/n}\mathbb{E}|S_1|} -1\biggr\rrvert >0.5
\biggr)
\\
&&\quad\qquad{} +\mathbb {P} \biggl(n^{1/\beta}\bigl|\Delta_i^n X\bigr|>0.5\alpha\sqrt{\frac{\pi}{2}}\sigma _{((j-1)k_n)/n}\mathbb{E}|S_1|n^{1/2-\varpi}
\biggr).
\end{eqnarray*}
From here we can use the bounds in (\ref{appecdf4}) and (\ref
{appecdf7}) to conclude
%
%
\begin{eqnarray}
\label{appecdf21}
\mathbb{P} \Bigl(\bigl|\Delta_i^n X\bigr|> \alpha
\sqrt{\widehat {V}{}^n_{j}}n^{-\varpi} \Bigr)
\leq\frac{K}{n^{\iota}}
\nonumber\\[-8pt]\\[-8pt]
\eqntext{\mbox{for some sufficiently small }\iota>0}
\end{eqnarray}
and hence the convergence in (\ref{appecdf19}) holds which implies
the result in (\ref{eqecdf1}).

\subsection{Proof of Theorem~\texorpdfstring{\protect\ref{teonoise}}{2}}
The proof follows the same steps as that of Theorem~\ref{teoecdf}. We
denote with $\eta_n$ a deterministic sequence depending only on $n$ and
vanishing as $n\rightarrow\infty$. Then, using the triangular
inequality and successive conditioning, we have
%
%
\begin{eqnarray}
\label{appnoise1} \mathbb{P} \biggl(\biggl\llvert \frac{1}{n}
\widehat{V}_j^n-\mu^2\biggr\rrvert +\biggl
\llvert \frac{1}{n}\widehat{V}_j^n(i)-
\mu^2\biggr\rrvert \geq\eta_n \biggr)&\leq& K
\frac{n^{-1/2}}{\eta_n},
\\
\label{appnoise2} \mathbb{P} \bigl((\varepsilon_{i/n}-\varepsilon_{(i-1)/n})
1_{ \{|\Delta_i^n X|> \alpha\sqrt{\widehat{V}{}^n_{j}}n^{-\varpi
} \}}\geq\eta_n \bigr)&\leq& K\frac{n^{\varpi-1/2}}{\eta
_n^{\iota}}.
\end{eqnarray}
We denote
%
%
\begin{eqnarray*}
\label{appnoise3}
\mathcal{B}_{i,j}^n
&=& \biggl\{\omega\dvtx  \bigl
\llvert \Delta_i^nX^*1_{ \{
|\Delta_i^n X|\leq\alpha\sqrt{\widehat{V}{}^n_{j}}n^{-\varpi}
\}}-(
\varepsilon_{i/n}-\varepsilon_{(i-1)/n})\bigr\rrvert >\eta_n
\\
&&\hspace*{146pt}{} \cup \biggl\llvert \frac{\sqrt{\widehat{V}_j^n(i)}}{{\sqrt{n}\mu}} -1\biggr\rrvert >\eta_n
\biggr\}
\end{eqnarray*}
for $i=(j-1)k_n+1,\ldots,(j-1)k_n+m_n$ and $j=1,\ldots,\lfloor
n/k_n\rfloor$. We set $\eta_n = n^{-x}$ for $0<x<\frac{1}{\iota
}(1/2-\varpi)\bigwedge1/2$. With this choice
%
%
\[
\label{appnoise4} \frac{1}{\lfloor n/k_n\rfloor m_n}\sum_{j=1}^{\lfloor n/k_n\rfloor}
\sum_{i=(j-1)k_n+1}^{(j-1)k_n+m_n}\mathbb{P} \bigl(\mathcal
{B}_{i,j}^n \bigr) = o(1).\vadjust{\goodbreak}
\]
Therefore, for any compact subset $\mathcal{A}$ of $(-\infty,0)$, we have
%
%
\[
\label{appnoise5} \sup_{\tau\in\mathcal{A}}\bigl|\widetilde{F}_n(\tau) -
\widehat {G}_n(\tau)\bigr| = o_p(1),
\]
where we denote
%
%
\begin{eqnarray*}
\label{appnoise6} \widehat{G}_n(\tau) &=& \frac{1}{\lfloor n/k_n\rfloor m_n}
\\[2pt]
&&{}\times \sum_{j=1}^{\lfloor n/k_n\rfloor} \sum_{i=(j-1)k_n+1}^{(j-1)k_n+m_n}1
\biggl\{ \frac{\sqrt{n}\Delta_i^nX}{\sqrt{\widehat{V}{}^n_{j}(i)}}1_{ \{
|\Delta
_i^n X|\leq\alpha\sqrt{\widehat{V}{}^n_{j}}n^{-\varpi} \}}\leq \tau \biggr\}1_{\{(\mathcal{B}_{i,j}^n)^c\}}.
\end{eqnarray*}
Taking into account the definition of the set $\mathcal{B}_{i,j}^n$,
we get
\[
\cases{ \displaystyle \widehat{G}_n(\tau)\geq
\frac{1}{\lfloor n/k_n\rfloor m_n}
\vspace*{7pt}
\cr
\hspace*{40pt}{}\times
\displaystyle\sum_{j=1}^{\lfloor n/k_n\rfloor}  \sum
_{i=(j-1)k_n+1}^{(j-1)k_n+m_n}1 \biggl\{\frac{1}{\mu} (
\varepsilon_{i/n}-\varepsilon_{(i-1)/n} )\leq\tau(1-\eta_n)-
\frac{\eta_n}{\mu} \biggr\},
\vspace*{7pt}
\cr
\displaystyle \widehat{G}_n(
\tau)\leq\frac{1}{\lfloor n/k_n\rfloor m_n}
\vspace*{7pt}
\cr
\hspace*{40pt}{}\times
\displaystyle\sum_{j=1}^{\lfloor n/k_n\rfloor}
\sum_{i=(j-1)k_n+1}^{(j-1)k_n+m_n}1 \biggl\{\frac{1}{\mu}
(\varepsilon_{i/n}-\varepsilon_{(i-1)/n} )\leq\tau(1+
\eta_n)+\frac{\eta_n}{\mu} \biggr\}.}
\]
From here we can proceed exactly in the same way as in the proof of
Theorem~\ref{teoecdf} to show that $\widehat{G}_n(\tau)\stackrel
{\mathbb{P}}{\longrightarrow} F_{\varepsilon}(\tau)$ locally uniformly in
$\tau$. Hence we need only show $\frac{N^n(\alpha,\varpi)}{\lfloor
n/k_n\rfloor m_n} \stackrel{\mathbb{P}}{\longrightarrow} 1$ as
$n\rightarrow\infty$. This follows from
%
%
\begin{eqnarray*}
\label{appnoise8}
&& \mathbb{P} \Bigl(\bigl|\Delta_i^nX^*\bigr|>\alpha
\sqrt{\widehat {V}{}^n_{j}}n^{-\varpi} \Bigr)
\\
&&\qquad \leq\mathbb{P} \biggl(\biggl\llvert \frac
{\sqrt {\widehat{V}{}^n_{j}}}{\sqrt{n}\mu}\biggr\rrvert >0.5 \biggr)+
\mathbb {P} \bigl(\bigl|\Delta_i^nX^*\bigr|>0.5\alpha\mu
n^{1/2-\varpi} \bigr)
\\
&&\qquad \leq \frac
{K}{n^{\iota}}\qquad\mbox{for some sufficiently small }\iota >0,
\end{eqnarray*}
which can be shown using (\ref{appnoise1}), the fact that the noise
term has a finite first moment and the Burkholder--Davis--Gundy
inequality.

\subsection{Proof of Theorem~\texorpdfstring{\protect\ref{teosv}}{3}} As in the proof of
Theorem~\ref{teoecdf}, without loss of generality we will assume $\tau<0$.
First, given the fact that $m_n/k_n\rightarrow0$, it is no\vadjust{\goodbreak} limitation
to assume $k_n-m_n>2$, and we will do so henceforth. Here we need to
make some additional decomposition of the difference $\widetilde{V}{}^n_j
-\widebar{V}_j^n$. It is given by the following:
%
%
\begin{eqnarray}
\widetilde{V}{}^n_j -\widebar{V}_j^n &=& R_j^{(1)}+R_j^{(2)}+R_j^{(3)}+R_j^{(4)},\qquad j=1,\ldots,\lfloor n/k_n\rfloor,\nonumber\hspace*{-23pt}
\\
R_j^{(1)} &=& \frac{n}{k_n-1}\frac{\pi}{2}\nonumber\hspace*{-23pt}
\\
&&{}\times \sum_{i=(j-1)k_n+2}^{jk_n} \bigl[\bigl(\bigl|\Delta_{i-1}^nB\bigr|\bigl|\Delta_i^nB\bigr|-\sigma_{(i-2)\Delta_n}^2\bigl|\Delta _{i-1}^nW\bigr|\bigl|\Delta_i^nW\bigr| \bigr)\nonumber\hspace*{-23pt}
\\
&&\hspace*{66pt}{} +(\sigma _{(i-2)\Delta_n} - \sigma_{((j-1)k_n)/n})^2\bigl|
\Delta_{i-1}^nW\bigr|\bigl|\Delta _i^nW\bigr|\bigr],\nonumber\hspace*{-23pt}
\\
R_j^{(2)} &=& 2\frac{n}{k_n-1}\frac{\pi}{2} \sigma_{((j-1)k_n)/n}\nonumber\hspace*{-23pt}
\\
&&{}\times \sum_{i=(j-1)k_n+2}^{jk_n}\biggl[ \sigma_{(i-2)\Delta_n} - \sigma_{((j-1)k_n)/n} \nonumber\hspace*{-23pt}
\\
&&\hspace*{66pt}{}- \int_{((j-1)k_n)/n}^{(i-2)/n}\tilde{\sigma}_{((j-1)k_n)/n}\,dW_u\nonumber\hspace*{-23pt}
\\
&&\hspace*{66pt}{}  - \int_{((j-1)k_n)/n}^{(i-2)/n}\tilde{\sigma}'_{((j-1)k_n)/n}\,dW'_u \biggr] \bigl|
\Delta_{i-1}^nW\bigr|\bigl|\Delta_i^nW\bigr|,\hspace*{-23pt}
\\
R_j^{(3)} &=& \frac{2}{k_n-1}\frac{\pi}{2} \sigma_{((j-1)k_n)/n}\nonumber\hspace*{-23pt}
\\
&&{}\times \sum_{i=(j-1)k_n+2}^{jk_n}\biggl[ \int_{((j-1)k_n)/n}^{(i-2)/n}\tilde{\sigma}_{((j-1)k_n)/n}\,dW_u \nonumber\hspace*{-23pt}
\\
&&\hspace*{66pt}{}+ \int_{((j-1)k_n)/n}^{(i-2)/n}\tilde{\sigma}'_{((j-1)k_n)/n}\,dW'_u\biggr]\nonumber\hspace*{-23pt}
\\
&&\hspace*{61pt}{}\times \biggl(n\bigl|\Delta _{i-1}^nW\bigr|\bigl|\Delta_i^nW\bigr|-\frac{2}{\pi} \biggr),\nonumber\hspace*{-23pt}
\\
R_j^{(4)} &=& \frac{2}{k_n-1}\sigma_{((j-1)k_n)/n}\nonumber\hspace*{-23pt}
\\
&&{}\times \sum_{i=(j-1)k_n+2}^{jk_n}\biggl[ \int_{((j-1)k_n)/n}^{(i-2)/n}\tilde{\sigma}_{((j-1)k_n)/n}\,dW_u\nonumber\hspace*{-23pt}
\\
&&\hspace*{66pt}{}  + \int_{((j-1)k_n)/n}^{(i-2)/n}\tilde{\sigma}'_{((j-1)k_n)/n}\,dW'_u\biggr].\nonumber\hspace*{-23pt}
\end{eqnarray}
For $i=(j-1)k_n+1,\ldots,jk_n-2$ we denote the component of $R_j^{(4)}$
that does not contain the increments $\Delta_i^nW$ and $\Delta_i^nW'$ with
%
%
\begin{eqnarray*}
\label{appsv6} \widetilde{R}_{i,j}^{(4)} &=&
R_j^{(4)}
\\
&&{}- \frac{2}{k_n-1} \sigma_{((j-1)k_n)/n}(jk_n-i-1)
\\
&&\quad{} \times\biggl[ \int_{(i-1)/n}^{i/n}\tilde{
\sigma}_{((j-1)k_n)/n}\,dW_u + \int_{(i-1)/n}^{i/n}
\tilde{\sigma}'_{((j-1)k_n)/n}\,dW'_u
\biggr].
\end{eqnarray*}
We decompose analogously the difference $\widetilde{V}{}^n_j(i)
-\widebar
{V}_j^n(i)$ into $R_j^{(k)}(i)$ for $k=1,\ldots,4$ and $\widetilde
{R}_{i,j}^{(4)}(i)$ is the component of $R_{j}^{(4)}(i)$ that does not
contain the increments $\Delta_i^nW$ and $\Delta_i^nW'$. We further
denote for $i=(j-1)k_n+1,\ldots,(j-1)k_n+m_n$ and $j = 1,\ldots,\lfloor
n/k_n \rfloor$,
%
%
\begin{eqnarray*}
\label{appsv7}
\xi_{j}^n(1) &=& \frac{\widehat{V}_j^n(i)-\sigma_{((j-1)k_n)/n}^2}{2\sigma^2_{((j-1)k_n)/n}},
\qquad
\xi_{j}^n(2) = \frac{ (\widehat{V}_j^n(i)-\sigma_{((j-1)k_n)/n}^2
)^2}{8\sigma_{((j-1)k_n)/n}^4},
\\[3pt]
\tilde{\xi}_{i,j}^n(1) &=& \frac{\widebar{V}_j^n(i)+\widetilde
{R}_{i,j}^{(4)}(i)-\sigma_{((j-1)k_n)/n}^2}{2\sigma_{((j-1)k_n)/n}^2},
\\[3pt]
\tilde{\xi}_{i,j}^n(2) &=& \frac{(\widebar
{V}_j^n(i)+\widetilde{R}_{i,j}^{(4)}(i)-\sigma_{((j-1)k_n)/n}^2 )^2}{8\sigma_{((j-1)k_n)/n}^4},
\\[3pt]
\bar{\xi}_{j}^n(1) &=& \frac{\widebar{V}_j^n+R_j^{(4)}-\sigma
_{((j-1)k_n)/n}^2}{2\sigma_{((j-1)k_n)/n}^2},
\qquad
\bar {\xi }_{j}^n(2) = \frac{ (\widebar{V}_j^n+R_j^{(4)}-\sigma_{((j-1)k_n)/n}^2 )^2}{8\sigma_{((j-1)k_n)/n}^4},
\\[3pt]
\hat {\xi}_{j}^n(1) &=& \frac{\widebar{V}_j^n-\sigma_{((j-1)k_n)/n}^2}{2\sigma_{((j-1)k_n)/n}^2},
\qquad
\hat{\xi }_{j}^n(2) = \frac{ (\widebar{V}_j^n-\sigma_{((j-1)k_n)/n}^2 )^2}{8\sigma_{((j-1)k_n)/n}^4},
\\[3pt]
\xi_{i,j}^n(3) &=& \frac{\sqrt{n}\Delta_i^nW}{\sigma_{((j-1)k_n)/n}} \bigl[\tilde{\sigma}_{((j-1)k_n)/n}(W_{(i-1)/n}-W_{((j-1)k_n)/n})
\\[3pt]
&&\hspace*{55pt}{} + \tilde{\sigma}'_{((j-1)k_n)/n}\bigl(W'_{(i-1)/n}-W'_{((j-1)k_n)/n}\bigr) \bigr],
\\[3pt]
\xi_{i,j}^n(4) &=& 1+\frac{1}{\sigma_{((j-1)k_n)/n}} \bigl[\tilde{\sigma}_{((j-1)k_n)/n}(W_{(i-1)/n}-W_{((j-1)k_n)/n})
\\[3pt]
&&\hspace*{73pt}{}+ \tilde{\sigma}'_{((j-1)k_n)/n}\bigl(W'_{(i-1)/n}-W'_{((j-1)k_n)/n}
\bigr) \bigr].
\end{eqnarray*}
With this notation we set for $i=(j-1)k_n+1,\ldots,(j-1)k_n+m_n$ and
$j=1,\ldots,\lfloor n/k_n\rfloor$
%
%
\begin{eqnarray*}
\label{appsv8} \chi_{i,j}^n(1) &=& -\sqrt{n}
\frac{1}{\sigma_{((j-1)k_n)/n}}
\\
&&\quad{}\times \biggl(\Delta_i^nA+
\Delta_i^nY+\int_{(i-1)/n}^{i/n}
(\sigma _{u}-\sigma_{(i-1)/n} )\,dW_u
\biggr)1_{ \{|\Delta_i^n
X|\leq\alpha\sqrt{\widehat{V}{}^n_{j}}n^{-\varpi} \}}
\\
&&{} +\bigl(\sqrt {n}\Delta_i^nW+\xi_{i,j}^n(3)
\bigr)1_{ \{|\Delta_i^n X|> \alpha
\sqrt
{\widehat{V}{}^n_{j}}n^{-\varpi} \}}
\\
&&{} - \biggl(\frac{\sqrt {n}\Delta
_i^nW}{\sigma_{((j-1)k_n)/n}}(\sigma_{(i-1)/n}-\sigma _{((j-1)k_n)/n}) -
\xi_{i,j}^n(3) \biggr)1_{ \{|\Delta_i^n X|\leq\alpha\sqrt
{\widehat{V}{}^n_{j}}n^{-\varpi} \}},
\nonumber
\\
\label{appsv8a} \chi_{i,j}^n(2) &=& \biggl(
\frac{\sqrt{\widehat{V}_j^n(i)}}{\sigma
_{((j-1)k_n)/n}}-1-\xi_{j}^n(1)+\xi_{j}^n(2)
\biggr)
\\
&&{} +\bigl(\xi_{j}^n(1) - \xi _{j}^n(2)
- \tilde{\xi}_{i,j}^n(1)+\tilde{\xi
}_{i,j}^n(2)\bigr).
\nonumber
\end{eqnarray*}
Finally, we denote
%
%
\begin{eqnarray*}
\label{appsv9} \widehat{G}_n(\tau) &=& \frac{1}{\lfloor n/k_n\rfloor m_n}
\\
&&{}\times \sum
_{j=1}^{\lfloor n/k_n\rfloor} \sum_{i=(j-1)k_n+1}^{(j-1)k_n+m_n}
1 \biggl(\sqrt{n}\frac{\Delta_i^nX}{\sigma_{((j-1)k_n)/n}} 1_{
\{
|\Delta_i^n X|\leq\alpha\sqrt{\widehat{V}{}^n_{j}}n^{-\varpi}
\}}
\\
&&\hspace*{103pt} \leq \tau\frac{\sqrt{\widehat{V}_j^n(i)}}{\sigma
_{((j-1)k_n)/n}}-\chi_{i,j}^n(1)-
\tau\chi_{i,j}^n(2) \biggr)
\\
& =& \frac
{1}{\lfloor n/k_n\rfloor m_n}
\\
&&{}\times \sum_{j=1}^{\lfloor n/k_n\rfloor} \sum
_{i=(j-1)k_n+1}^{(j-1)k_n+m_n}1 \bigl(\sqrt{n}
\Delta_i^nW \leq\tau +\tau \tilde{\xi}_{i,j}^n(1)-\tau\tilde{\xi}_{i,j}^n(2)
-\xi _{i,j}^n(3) \bigr).
\nonumber
\end{eqnarray*}
The proof consists of three parts: the first is showing the
negligibility of $k_n(\widetilde{F}_n(\tau) - \widehat{G}_n(\tau))$,
the second is deriving the limiting behavior of $\widehat{G}_n(\tau) -
\Phi(\tau)$ and third part is showing negligibility of $k_n(\widehat
{F}_n(\tau) - \widetilde{F}_n(\tau))$.

\subsubsection{The difference \texorpdfstring{$\widetilde{F}_n(\tau)-\widehat{G}_n(\tau)$}
{Fn(tau) - Gn(tau)}}
We first collect some preliminary results that we then make use of in
analyzing $\widetilde{F}_n(\tau) - \widehat{G}_n(\tau)$. We start with
$\max_{i=1,\ldots,n}|\Delta_i^n B|$. Using maximal inequality we have
%
%
\begin{equation}
\label{appsv14} \mathbb{E}\Bigl(\max_{i=1,\ldots,n}\bigl|
\Delta_i^n B\bigr|^p\Bigr) \leq
Kn^{1-p/2}\qquad \forall p>0.
\end{equation}
Next, using Assumption \hyperref[assSB]{SB} (in particular that jumps are of finite
activity), we have
%
%
\begin{eqnarray}
\label{appsv15}
\mathbb{P} \biggl( \int_{((j-1)k_n)/n}^{(jk_n)/n}
\int_{E}1 \bigl(\delta^{\phi}(z,x)\neq0 \bigr)
\mu(dz,dx)\geq1 \biggr)\leq K\frac{k_n}{n},
\nonumber\\[-8pt]\\[-8pt]
\eqntext{\phi= Y, \sigma, \tilde{\sigma}\mbox{ and }\tilde{\sigma}'.}
\end{eqnarray}
We now provide bounds for the elements of $\chi_{i,j}^n(1)$ and $\chi
_{i,j}^n(2)$. In what follows we denote with $\eta_n$ some
deterministic sequence of positive numbers that depends only on $n$. We
first have (recall the definition of $\bar{\sigma}_t$)
%
%
\begin{eqnarray*}
\label{appsv16}
&& \mathbb{P} \biggl( \sqrt{n}\biggl\llvert \int_{(i-1)/n}^{i/n}(
\sigma _u-\sigma_{(i-1)/n})\,dW_u\biggr\rrvert
\geq\eta_n \biggr)
\\
&&\qquad \leq \mathbb {P} \biggl( \int_{((j-1)k_n)/n}^{(jk_n)/n}
\int_{E}1 \bigl(\delta^Y(s,x)\neq0 \bigr)
\mu(ds,dx)\geq1 \biggr)
\\
&&\quad\qquad{} +\mathbb {P} \biggl( \sqrt{n}\biggl\llvert \int_{(i-1)/n}^{i/n}(
\bar{\sigma }_u-\bar{\sigma}_{(i-1)/n})\,dW_u
\biggr\rrvert \geq\eta _n \biggr).
\nonumber
\end{eqnarray*}
For the second term on the right-hand side of the above inequality, we
can use Chebyshev's inequality as well as Burkholder--Davis--Gundy
inequality, to get for~$\forall p\geq2$
%
%
\begin{eqnarray*}
\label{appsv17}
&& \mathbb{P} \biggl( \sqrt{n}\biggl\llvert \int_{(i-1)/n}^{i/n}(
\bar{\sigma}_u-\bar{\sigma}_{(i-1)/n})\,dW_u
\biggr\rrvert \geq\eta_n \biggr)
\\
&&\qquad \leq \frac{n^{p/2}\mathbb{E}\llvert \int_{(i-1)\Delta_n}^{i\Delta
_n}(\bar{\sigma}_u-\bar{\sigma}_{(i-1)/n})^2\,du\rrvert ^{p/2}}{\eta
_n^{p}}.
\end{eqnarray*}
Therefore, applying again the Burkholder--Davis--Gundy inequality, we have
altogether
%
\begin{eqnarray}
\label{appsv18}
&& \mathbb{P} \biggl( \sqrt{n}\biggl\llvert \int_{(i-1)/n}^{i/n}(
\sigma _u-\sigma_{(i-1)/n})\,dW_u\biggr\rrvert
\geq\eta_n \biggr)
\nonumber\\[-8pt]\\[-8pt]
&&\qquad \leq K \biggl[ \biggl(\frac{k_n}{n} \biggr)
\vee \biggl(\frac{1}{n^{p/2}\eta_n^{p}} \biggr) \biggr]\qquad \forall p>0.\nonumber
\end{eqnarray}
Similar calculations (using the fact that $\tilde{\sigma}_t$ and
$\tilde{\sigma}'_t$ are It\^o semimartingales), yields for
$\forall p>0$
%
%
\begin{eqnarray}\label{appsv18a}
&& \mathbb{P} \biggl( \biggl\llvert \frac{\sqrt{n}\Delta_i^nW}{\sigma_{((j-1)k_n)/n}} (
\sigma_{(i-1)/n}-\sigma_{((j-1)k_n)/n} )-\xi_{i,j}^n(3)
\biggr\rrvert \geq\eta_n \biggr)
\nonumber\\[-8pt]\\[-8pt]
&&\qquad \leq K \biggl[ \biggl(
\frac{k_n}{n} \biggr)\vee \biggl(\frac{k_n}{n\eta_n}
\biggr)^p \biggr].\nonumber
\end{eqnarray}
Next, applying Chebyshev's inequality and the elementary
$|\sum_i|a_i||^p\leq\break \sum_i|a_i|^p$ for $p\in(0,1]$, we get
%
%
\begin{eqnarray}
\label{appsv19} \mathbb{P} \bigl(\sqrt{n}\bigl|\Delta_i^nY\bigr|
\geq\eta_n \bigr)&\leq& \frac
{n^{\iota/2}\mathbb{E} (\int_{(i-1)\Delta_n}^{i\Delta_n}\int_{E}|\delta^Y(s,x)|\mu(ds,dx)  )^{\iota}}{\eta_n^{\iota}}\nonumber
\\
&\leq&\frac{n^{\iota/2}\mathbb{E} (\int_{(i-1)\Delta
_n}^{i\Delta
_n}\int_{E}|\delta^Y(s,x)|^{\iota}\mu(ds,dx)  )}{\eta
_n^{\iota}}
\\
&\leq& Kn^{-1+\iota/2}\eta_n^{-\iota}\qquad \forall\iota
\in(0,1).\nonumber
\end{eqnarray}
Further, Chebyshev's inequality and the boundedness of $a_t$ easily implies
%
%
\begin{equation}
\label{appsv20} \mathbb{P} \bigl(\sqrt{n}\bigl|\Delta_i^nA\bigr|
\geq\eta_n \bigr)\leq \frac
{n^{p/2}\mathbb{E}(|\Delta_i^nA|^p)}{\eta^p_n}\leq K\frac
{1}{n^{p/2}\eta^p_n}.
\end{equation}
We turn next to the difference $\widehat{V}_j^n-\dot{V}_j^n$. Using the
triangular inequality and successive conditioning, we have
%
%
\begin{eqnarray*}
\label{appsv21} && \mathbb{P} \bigl(\bigl|\widehat{V}_j^n-
\dot{V}_j^n\bigr|\geq\eta_n \bigr)
\\
&&\qquad \leq\mathbb{P} \biggl( 2\frac{n}{k_n}\frac{\pi}{2}\max
_{i=1,\ldots,n}\bigl|\Delta_i^nA+
\Delta_i^nB\bigr|
\\
&&\hspace*{42pt}{}\times \int_{((j-1)k_n)/n}^{(jk_n)/n}
\int_{E}\bigl(\bigl|\delta^Y(s,x)\bigr|\vee1\bigr)
\mu(ds,dx)\geq\frac{\eta
_n}{2} \biggr)+\frac{K}{n\eta_n}.
\nonumber
\end{eqnarray*}
From here we have
%
%
\begin{eqnarray*}
\label{appsv23} &&\mathbb{P} \biggl( 2\frac{n}{k_n}\frac{\pi}{2}\max
_{i=1,\ldots,n}\bigl|\Delta _i^nA+
\Delta_i^nB\bigr|\int_{((j-1)k_n)/n}^{(jk_n)/n}
\int_{E}\bigl(\bigl|\delta^Y(s,x)\bigr|\vee1\bigr)
\mu(ds,dx)\geq\frac{\eta_n}{2} \biggr)
\\
&&\qquad\leq\mathbb{P} \biggl(\int_{((j-1)k_n)/n}^{(jk_n)/n}
\int_{E}\mu(ds,dx)\geq1 \biggr)\leq K\frac{k_n}{n}.
\nonumber
\end{eqnarray*}
Thus altogether we get
%
%
\begin{equation}
\label{appsv24} \mathbb{P} \bigl(\bigl|\widehat{V}_j^n-
\dot{V}_j^n\bigr|\geq\eta_n \bigr)\leq K \biggl(
\frac{1}{n\eta_n}\vee\frac{k_n}{n} \biggr).
\end{equation}

We continue next with the difference $\dot{V}_j^n - \widetilde{V}_j^n$.
Application of triangular inequality gives
%
%
\begin{eqnarray*}
\label{appsv25} &&\bigl|\Delta_{i-1}^nA+\Delta_{i-1}^nB\bigr|\bigl|
\Delta_i^nA+\Delta_i^nB\bigr|-\bigl|
\Delta _{i-1}^nB\bigr|\bigl|\Delta_i^nB\bigr|
\\
&&\qquad \leq \bigl|\Delta_{i-1}^nA+
\Delta_{i-1}^nB\bigr|\bigl|\Delta_i^nA\bigr|+\bigl|
\Delta _{i-1}^nA\bigr|\bigl|\Delta _i^nB\bigr|.
\nonumber
\end{eqnarray*}
Using this inequality and applying Chebyshev's inequality, we get
%
%
\begin{equation}
\label{appsv26} \mathbb{P} \bigl(\bigl|\dot{V}_j^n -
\widetilde{V}_j^n\bigr|\geq\eta_n \bigr)\leq K
\biggl(\frac{1}{\sqrt{n}\eta_n} \biggr)^p\qquad \forall p\geq1
\end{equation}
and this inequality can be further strengthened but suffices for our
analysis.\vadjust{\goodbreak}

Turning next to $R_j^{(1)}$, using the triangular inequality, the
Burkholder--Davis--Gundy inequality as well as (\ref{appsv15}), we
can easily get
%
%
\begin{eqnarray*}
\label{appsv27}
&&
\mathbb{P} \bigl(\bigl|R_j^{(1)}\bigr|\geq
\eta_n \bigr)
\\
&&\qquad \leq \mathbb {P} \biggl(\bigl|R_j^{(1)}\bigr|
\geq\eta_n, \int_{((j-1)k_n)/n}^{(jk_n)/n}\int
_{E}1\bigl(\delta^{\sigma}(s,x)\neq0\bigr)\mu(ds,dx)
\geq1 \biggr)
\\
&&\quad\qquad{} +\mathbb {P} \biggl(\bigl|R_j^{(1)}\bigr|\geq
\eta_n, \int_{((j-1)k_n)/n}^{(jk_n)/n}\int
_{E}1\bigl(\delta^{\sigma}(s,x)\neq0\bigr)\mu(ds,dx)=
0 \biggr)
\\
&&\qquad \leq K\frac{k_n}{n}+K \biggl(\frac{1}{\sqrt{n}\eta_n} \biggr)^p + K
\biggl(\frac{k_n}{n\eta_n} \biggr)^p\qquad \forall p\geq1.
\end{eqnarray*}
Similar calculations, and utilizing the fact that $\tilde{\sigma}_t$ $\tilde{\sigma}'_t$ are themselves It\^o semimartingales, yield
%
%
\begin{equation}
\label{appsv28} \mathbb{P} \bigl(\bigl|R_j^{(2)}\bigr|\geq
\eta_n \bigr)\leq K\frac{k_n}{n} + K \biggl(\frac{k_n}{n\eta_n}
\biggr)^p\qquad \forall p\geq1.
\end{equation}
Next, by splitting
%
%
\begin{eqnarray*}
\label{appsv29}
&& n\bigl|\Delta_{i-1}^nW\bigr|\bigl|\Delta_i^nW\bigr|-
\frac{2}{\pi}
\\
&&\qquad = \bigl|\sqrt{n}\Delta _{i-1}^nW\bigr| \biggl(\bigl|
\sqrt{n}\Delta_i^nW\bigr|-\sqrt{\frac{2}{\pi}} \biggr)+
\sqrt {\frac{2}{\pi}} \biggl(\bigl|\sqrt{n}\Delta_{i-1}^nW\bigr| -
\sqrt{\frac
{2}{\pi
}} \biggr),
\end{eqnarray*}
we can decompose $R_j^{(3)}$ into two discrete martingales. Then
applying the Burkholder--Davis--Gundy inequality, we get
%
%
\begin{equation}
\label{appsv30} \mathbb{P} \bigl(\bigl|R_j^{(3)}\bigr|\geq
\eta_n \bigr)\leq K \biggl(\frac
{1}{\sqrt {n}\eta_n} \biggr)^{p}\qquad \forall p\geq2.
\end{equation}
Next, we trivially have
%
%
\begin{equation}
\label{appsv31} \cases{ \displaystyle \mathbb{P} \bigl(\bigl|\widebar{V}_j^{n}-
\sigma_{((j-1)k_n)/n}^2\bigr|\geq\eta_n \bigr)\leq K \biggl(
\frac{1}{k_n\eta_n^2} \biggr)^p, \vspace*{6pt}
\cr
\displaystyle
\mathbb{P} \bigl(\bigl|\widetilde{V}_j^n-\widebar
{V}_j^{n}\bigr|\geq0.5\sigma_{((j-1)k_n)/n}^2
\bigr)\leq K\frac
{k_n}{n}, &\quad $\forall p\geq2$.}
\end{equation}
Further, application of the Burkholder--Davis--Gundy inequality gives
%
%
\begin{equation}
\label{appsv32} \qquad \cases{ \displaystyle \mathbb{E}\bigl|\widebar{V}_j^{n}-
\sigma_{((j-1)k_n)/n}^2\bigr|^p\leq\frac{K}{k_n^{p/2}},
\vspace*{5pt}
\cr
\displaystyle \mathbb{E}_{(j-1)k_n}^n
\bigl(R_j^{(4)}\bigr)=0,\qquad \mathbb{E}\bigl|R_j^{(4)}\bigr|^p
\leq K \biggl(\frac
{k_n}{n} \biggr)^{p/2}, &\quad $\forall p
\geq2$.}
\end{equation}
The results in (\ref{appsv24})--(\ref{appsv32}) continue to hold
when $\widehat{V}_j^{n}$, $\dot{V}_j^n$, $\widetilde{V}_j^n$,
$\widebar{V}_j^n$, $R_j^{(1)}$, $R_j^{(2)}$ and $R_j^{(3)}$ are replaced with
$\widehat{V}_j^{n}(i)$, $\dot{V}_j^n(i)$, $\widetilde{V}_j^n(i)$,
$\widebar{V}_j^n(i)$, $R_j^{(1)}(i)$, $R_j^{(2)}(i)$ and
$R_j^{(3)}(i)$, respectively.

Further, using the Burkholder--Davis--Gundy inequality for discrete
martingales [note that $\widebar{V}_j^{n} - \widebar{V}_j^{n}(i)$ can
be decomposed into discrete martingales and terms whose $p$th moment is
bounded by $K/k_n^p$], we have
%
%
\begin{eqnarray}
\label{appsv33} \qquad\mathbb{E} \bigl(\bigl|R_j^{(4)} - \widetilde
{R}_{i,j}^{(4)}\bigr|^p+\bigl|\widetilde
{R}_{i,j}^{(4)} - \widetilde{R}_{i,j}^{(4)}(i)\bigr|^p
\bigr)&\leq& K \biggl(\frac{1}{\sqrt{n}} \biggr)^p\qquad \forall p>0,
\\
\bigl\llvert \mathbb{E}_{(j-1)k_n}^n \bigl(
\widebar{V}_j^{n} - \widebar {V}_j^{n}(i)
\bigr)\bigr\rrvert &\leq&\frac{K}{k_n^2},
\nonumber\\[-8pt]\label{appsv34}  \\[-8pt]
\mathbb {E}\bigl|\widebar {V}_j^{n} - \widebar{V}_j^{n}(i)\bigr|^p
&\leq& K \biggl(\frac{1}{k_n} \biggr)^p\qquad \forall p\geq1.\nonumber
\end{eqnarray}
Now we can use the above results for the components of $\widehat
{V}_j^n(i) -\sigma_{((j-1)k_n)/n}^2$, to analyze the first term
in $\chi_{i,j}^n(2)$ involving $\sqrt{\widehat{V}_j^n(i)}-\sigma
_{((j-1)k_n)/n}$. We make use of the following algebraic inequality:
%
%
\[
\label{appsv35} \biggl\llvert \sqrt{x}-\sqrt{y}-\frac{x-y}{2\sqrt{y}}+
\frac
{(x-y)^2}{8y\sqrt {y}}\biggr\rrvert \leq\frac{(x-y)^4}{8y^{7/2}} + \frac
{|x-y|^3}{2y^{5/2}}\nonumber
\]
for every $x\geq0$ and $y>0$. Using this inequality with $x$ and $y$
replaced with $\widehat{V}_j^n(i)$ and $\sigma_{((j-1)k_n)/n}^2$,
respectively, as well the bounds in (\ref{appsv24})--(\ref
{appsv34}), we get
%
%
\begin{eqnarray}
\label{appsv36} &&\mathbb{P} \biggl( \biggl\llvert \frac{\sqrt{\widehat{V}_j^n(i)}}{\sigma
_{((j-1)k_n)/n}}-1-
\xi_{j}^n(1)+\xi_{j}^n(2)\biggr
\rrvert \geq\eta_n \biggr)
\nonumber\\[-8pt]\\[-8pt]
&&\qquad\leq K \biggl[ \frac{1}{n\eta_n^{1/3}} \vee\frac{k_n}{n}
\vee\frac
{1}{\eta_n^{p/3}[n^{p/2}\wedge(n/k_n)^{p/2}]} \vee\frac{1}{\eta
_n^{2p/3}k_n^p} \biggr]\nonumber
\end{eqnarray}
for $\forall p\geq1$ and $\forall\iota>0$. Similarly, using the
following inequality:
%
%
\[
\label{appsv37} \mathbb{P} \bigl(\bigl|x^2-y^2\bigr|\geq\varepsilon
\bigr)\leq\mathbb{P} \bigl(|x-y|^2\geq0.5\varepsilon \bigr)+\mathbb{P}
\bigl(2|y|\geq K \bigr)+\mathbb{P} \bigl(|x-y|\geq0.5\varepsilon/K \bigr)
\]
for any random variables $x$ and $y$ and constants $\varepsilon>0$ and
$K>0$, together with the bounds in (\ref{appsv24})--(\ref
{appsv34}), we have
%
%
\begin{eqnarray}
\label{appsv38} &&\mathbb{P} \bigl( \bigl|\xi_{j}^n(1) -
\xi_{j}^n(2) - \tilde{\xi }_{i,j}^n(1)+
\tilde{\xi}_{i,j}^n(2)\bigr|\geq\eta_n \bigr)
\nonumber\\[-8pt]\\[-8pt]
&&\qquad\leq K \biggl[ \frac{1}{n\eta_n} \vee
\frac
{k_n}{n}\vee\frac{1}{\eta_n^{p}[n^{p/2}\wedge(n/k_n)^p]} \biggr]\nonumber
\end{eqnarray}
for every $p\geq1$ and arbitrary small $\iota>0$.

We finally provide a bound for the second term in $\chi_{i,j}^n(1)$. We
can use Chebyshev inequality as well as H\"{o}lder's inequality to get
%
%
\begin{eqnarray}
\label{appsv39} \qquad &&\mathbb{P} \Bigl(\bigl(\sqrt{n}\sigma_{((j-1)k_n)/n}\bigl|\Delta
_i^nW\bigr|+\bigl|\xi _{i,j}^n(3)\bigr|\bigr)1
\Bigl(\bigl|\Delta_i^nX\bigr|> \alpha\sqrt{\widehat
{V}_j^n}n^{-\varpi} \Bigr)\geq\eta_n
\Bigr)
\nonumber\\[-8pt]\\[-8pt]
&&\qquad         \leq K\frac{ [\mathbb{P}
(|\Delta_i^nX|> \alpha\sqrt{\widehat{V}_j^n}n^{-\varpi}
)]^{1/(1+\iota)}}{\eta_n^{\iota}}.\nonumber
\end{eqnarray}
We can further write
%
%
\begin{eqnarray*}
\label{appsv40}
\mathbb{P} \Bigl(\bigl|\Delta_i^nX\bigr|> \alpha
\sqrt{\widehat {V}_j^n}n^{-\varpi
} \Bigr)
&\leq& \mathbb{P} \Bigl(\bigl|\sqrt{\widehat{V}_j^n}-
\sigma _{((j-1)k_n)/n}\bigr|\geq0.5 \sigma_{((j-1)k_n)/n} \Bigr)
\\
&&{}
+ \mathbb{P} \bigl(\bigl|\Delta_i^nX\bigr|> 0.5\alpha
\sigma_{((j-1)k_n)/n}n^{-\varpi} \bigr).
\end{eqnarray*}
From here we can use the bounds in (\ref{appsv24})--(\ref{appsv34})
as well as (\ref{appsv19}) and conclude
%
%
\begin{eqnarray}
\label{appsv41} &&\mathbb{P} \Bigl(\sqrt{n}\sigma_{((j-1)k_n)/n}\bigl|\Delta
_i^nW\bigr|1 \Bigl(\bigl|\Delta_i^nX\bigr|>
\alpha\sqrt{\widehat{V}_j^n}n^{-\varpi}
\Bigr)\geq\eta _n \Bigr)
\nonumber\\[-8pt]\\[-8pt]
&&\qquad              \leq K \biggl(\frac{k_n}{n}
\biggr)^{1/(1+\iota)}\frac{1}{\eta_n^{\iota
}}\qquad \forall\iota>0.\nonumber
\end{eqnarray}

Combining the results in (\ref{appsv14}), (\ref{appsv18}), (\ref{appsv18a}), (\ref{appsv19}), (\ref{appsv20}), (\ref{appsv36}),
(\ref{appsv38}) and (\ref{appsv41}), we get
%
%
\begin{eqnarray*}
\label{appsv42} &&\mathbb{P} \bigl(\bigl(\bigl|\chi_{i,j}^n(1)\bigr|+\bigl|
\chi_{i,j}^n(2)\bigr|\bigr)>\eta_n \bigr)
\\
&&\qquad  \leq K \biggl[ \frac{1}{n\eta_n} \vee\frac{1}{\eta
_n^{p}[n^{p/2}\wedge(n/k_n)^p\wedge k_n^{3p/2}]}
\vee \biggl(\frac
{k_n}{n} \biggr)^{1/(1+\iota)}
\frac{1}{\eta_n^{\iota}} \biggr].
\nonumber
\end{eqnarray*}
From here, using the fact that the probability density of a standard
normal variable is uniformly bounded, we get
%
%
\begin{eqnarray*}
\label{appsv43} \hspace*{-6pt}&&\mathbb{E} \biggl|1 \Bigl(\sqrt{n}\Delta_i^nX
1_{ \{|\Delta_i^n
X|\leq\alpha\sqrt{\widehat{V}{}^n_{j}}n^{-\varpi} \}}\leq\tau \sqrt{\widehat{V}_j^n(i)}
\Bigr)
\\
\hspace*{-6pt}&&\hspace*{6pt}{}   - 1 \biggl(\sqrt{n}\frac{\Delta_i^nX}{\sigma_{((j-1)k_n)/n}} 1_{ \{|\Delta_i^n X|\leq\alpha\sqrt{\widehat
{V}{}^n_{j}}n^{-\varpi
} \}}\leq\tau\frac{\sqrt{\widehat{V}_j^n(i)}}{\sigma_{((j-1)k_n)/n}}-\chi_{i,j}^n(1)-\tau
\chi_{i,j}^n(2) \biggr) \biggr|
\\
\hspace*{-6pt}&&\hspace*{-4pt}\qquad \leq \mathbb{P} \bigl(\bigl(\bigl|\chi_{i,j}^n(1)\bigr|+\bigl|
\chi_{i,j}^n(2)\bigr|\bigr)>\eta_n \bigr)
\\
\hspace*{-6pt}&&\hspace*{-4pt}\qquad\quad{}+\mathbb{E}\biggl\llvert \Phi \biggl(\frac{\tau+\tau\tilde{\xi
}_{i,j}^n(1)-\tau\tilde{\xi}_{i,j}^n(2)+\eta_n(1+|\tau|)}{\xi
_{i,j}^n(4)} \biggr)
\\
\hspace*{-6pt}&&\hspace*{51pt}{} - \Phi\biggl(\frac{\tau+\tau\tilde{\xi
}_{i,j}^n(1)-\tau\tilde{\xi}_{i,j}^n(2)-\eta_n(1+|\tau|)}{\xi
_{i,j}^n(4)} \biggr) \biggr\rrvert
\\
\hspace*{-6pt}&&\hspace*{-4pt}\qquad \leq K\mathbb{P} \bigl(\bigl(\bigl|\chi _{i,j}^n(1)\bigr|+\bigl|
\chi_{i,j}^n(2)\bigr|\bigr)>\eta_n \bigr) + K
\eta_n|\tau|.
\nonumber
\end{eqnarray*}
Therefore, upon picking $\eta_n\propto n^{-q-\iota}$ for $\iota\in(0,
1/2-q)$ sufficiently small, we get finally for any compact subset
$\mathcal{A}$ of $(-\infty,0)$
%
%
\begin{equation}
\label{appsv44} \sup_{\tau\in\mathcal{A}}\bigl|\widetilde{F}_n(\tau) -
\widehat {G}_n(\tau)\bigr| = o_p \biggl(\frac{1}{k_n}
\biggr).
\end{equation}

\subsubsection{The asymptotic behavior of \texorpdfstring{$\widehat{G}_n(\tau)-\Phi(\tau)$}
{Gn(tau) - Phi(tau)}}
We have
%
%
\begin{eqnarray}
 \widehat{G}_n(\tau) - \Phi(\tau)&=&\sum
_{i=1}^5A_i^n, \nonumber
\\[-1pt]
A_1^n &=& \frac
{1}{\lfloor n/k_n\rfloor m_n}\sum_{j=1}^{\lfloor n/k_n\rfloor} \sum
_{i=(j-1)k_n+1}^{(j-1)k_n+m_n} \bigl[1 \bigl(\sqrt{n}
\Delta_i^nW\leq \tau \bigr)-\Phi(\tau) \bigr],\nonumber
\\[-1pt]
A_2^n &=& \frac{1}{\lfloor n/k_n\rfloor}\sum
_{j=1}^{\lfloor n/k_n\rfloor } \bigl(\Phi \bigl(\tau+\tau\bar{\xi}_{j}^n(1)-\tau\bar {\xi }_{j}^n(2)
\bigr) - \Phi(\tau) \bigr),\nonumber
\\[-1pt]
A_3^n &=& \frac{1}{\lfloor
n/k_n\rfloor m_n} \sum_{j=1}^{\lfloor n/k_n\rfloor} \sum
_{i=(j-1)k_n+1}^{(j-1)k_n+m_n}a_i^n,\nonumber
\\[-1pt]
a_i^n &=& 1 \biggl(\sqrt{n}
\Delta_i^nW \leq\frac{\tau+\tau\tilde
{\xi}_{i,j}^n(1)-\tau\tilde{\xi}_{i,j}^n(2)}{\xi_{i,j}^n(4)} \biggr) - 1 \bigl(
\sqrt{n}\Delta_i^nW\leq\tau \bigr)
\nonumber\\[-9pt]\\[-9pt]
&&{} +\Phi(\tau)- \Phi\biggl(\frac
{\tau+\tau\tilde{\xi}_{i,j}^n(1)-\tau\tilde{\xi
}_{i,j}^n(2)}{\xi_{i,j}^n(4)} \biggr),\nonumber
\\[-1pt]
A_4^n &=& \frac{1}{\lfloor n/k_n\rfloor m_n}
\sum
_{j=1}^{\lfloor
n/k_n\rfloor} \sum_{i=(j-1)k_n+1}^{(j-1)k_n+m_n}
\bigl[ \Phi \bigl(\tau +\tau\tilde{\xi}_{i,j}^n(1)-\tau
\tilde{\xi }_{i,j}^n(2) \bigr) \nonumber
\\[-1pt]
&&\hspace*{127pt}{}- \Phi \bigl(\tau+\tau
\bar{\xi}_{j}^n(1)-\tau\bar{\xi
}_{j}^n(2) \bigr) \bigr],\nonumber
\\[-1pt]
A_5^n &=& \frac{1}{\lfloor n/k_n\rfloor m_n}\sum
_{j=1}^{\lfloor
n/k_n\rfloor} \sum_{i=(j-1)k_n+1}^{(j-1)k_n+m_n}
\biggl[ \Phi \biggl(\frac
{\tau+\tau\tilde{\xi}_{i,j}^n(1)-\tau\tilde{\xi
}_{i,j}^n(2)}{\xi_{i,j}^n(4)} \biggr)\nonumber
\\[-1pt]
&&\hspace*{129pt}{} - \Phi \bigl(\tau+\tau
\tilde{\xi }_{i,j}^n(1)-\tau\tilde{\xi}_{i,j}^n(2) \bigr) \biggr].\nonumber
\end{eqnarray}
We first derive a bound for the order of magnitude of $A_3^n$, $A_4^n$
and $A_5^n$ and then analyze the limiting behavior of $A_1^n$ and $A_2^n$.
Using the independence of $\Delta_i^nW$, $\Delta_h^nW$, $\Delta_i^nW'$,
$\Delta_h^nW'$ from each\vspace*{1pt} other (for $i\neq h$) and $\mathcal
{F}_{((j-1)k_n)/n}$, the fact that $\xi_{i,j}^n(4)$ is adapted to
$\mathcal
{F}_{i-1}^n$ as well as successive conditioning, we have $\mathbb
{E} (a_i^na_h^n ) = 0$ for $|i-h|>k_n$. For $0<i-h\leq k_n$,
we can first split $a_h^n$ into a component in which the summand
including the $i$th increment $\Delta_i^nW$ is removed from
$\tilde{\xi}_{h,j}^n(1)$ and $\tilde{\xi}_{h,j}^n(2)$. We denote this part
of $a_h^n$ with $\bar{a}_h^n$ and the residual with $\tilde{a}_h^n = a_h^n-\bar{a}_h^n$. We\vadjust{\goodbreak} further denote with $\tilde{\xi}_{h,j}^{i,n}(1)$ and $\tilde{\xi}_{h,j}^{i,n}(2)$ the terms
$\tilde{\xi}_{h,j}^n(1)$ and $\tilde{\xi}_{h,j}^n(2)$ in which
the summand corresponding to $\Delta_i^nW$ is removed. Then using
successive conditioning, we have for $(j-1)k_n+1\leq h<i\leq(j-1)k_n+m_n$
%
%
\[
\label{appsv49} \mathbb{E} \bigl(a_i^n
\bar{a}_h^n \bigr) = 0,\qquad \mathbb{E}
\bigl(a_i^n \bigr)^2\leq K|\tau| \biggl(
\sqrt{\frac{k_n}{n}}\vee\frac
{1}{\sqrt{k_n}} \biggr).
\]
Further, we can use the triangular inequality for $a_i^n$ and
$\tilde{a}_h^n$, the bounds in (\ref{appsv32})--(\ref
{appsv34}), and get for $n$ sufficiently high
%
%
\begin{eqnarray*}
\label{appsv49a}
\mathbb{E}\bigl|a_i^n\tilde{a}_h^n\bigr|
&\leq& \mathbb{P} \biggl(\bigl|\tilde{\xi}_{i,j}^n(4)-1\bigr|> \biggl(
\frac
{k_n}{n} \biggr)^{1/2-\iota}\cup \bigl|\tilde{\xi}_{h,j}^n(4)-1\bigr|>
\biggl(\frac{k_n}{n} \biggr)^{1/2-\iota} \biggr)
\\[-1pt]
&&{} +\mathbb{P}\bigl(\bigl|\tilde{\xi}_{h,j}^n(1)-\tilde{\xi}_{h,j}^n(2)\bigr|>k_n^{-1/2+\iota}\bigr)
\\[-1pt]
&&{} +\mathbb{P} \bigl(\bigl|\tilde{\xi}_{i,j}^n(1)-
\tilde{\xi }_{i,j}^n(2)\bigr|>k_n^{-1/2+\iota}
\bigr)
\\[-1pt]
&&{} +\mathbb{P} \bigl(\bigl|\tilde {\xi }_{h,j}^n(1) -
\tilde{\xi}_{h,j}^{i,n}(1)-\tilde{\xi
}_{h,j}^n(2)+\tilde{\xi}_{h,j}^{i,n}(2)\bigr|>k_n^{-1+\iota}
\bigr)
\\[-1pt]
&&{} + \mathbb{P} \bigl(\sqrt{n}\Delta_i^nW\in2\tau
\bigl(1-k_n^{-1/2+\iota
}, 1+k_n^{-1/2+\iota}\bigr)
\\[-1pt]
&&\hspace*{28pt}{} \cap \sqrt{n}\Delta_h^nW\in2\tau \bigl(1-k_n^{-1+\iota
},1+k_n^{-1+\iota}
\bigr) \bigr)
+K\frac{|\tau|\vee\tau^2}{k_n}.
\end{eqnarray*}
Therefore, using again (\ref{appsv32})--(\ref{appsv34}), we have
%
%
\begin{equation}
\label{appsv50} A_3^n \leq K \bigl(\sqrt{|\tau|}\vee
\tau^2\bigr)\times \biggl(\frac{1}{\sqrt {\lfloor n/k_n\rfloor m_n}} \biggl(
\frac{1}{k_n} \biggr)^{1/4}\vee \frac
{1}{\sqrt{n}}
\biggr).
\end{equation}
For $A_4^n$, using a second-order Taylor expansion, the bounds in (\ref
{appsv31}), (\ref{appsv32}) and (\ref{appsv34}), as well as the
uniform boundedness of the probability density of the standard normal
distribution and its derivative, we get
%
%
\begin{equation}
\label{appsv51} \mathbb{E}\bigl|A_4^n\bigr|\leq K\bigl(|\tau|\vee
\tau^2\bigr) \biggl(\frac
{1}{k_n^{3/2}}\vee
\frac{1}{\sqrt{\lfloor n/k_n\rfloor}k_n} \biggr).
\end{equation}
Next, for $A_5^n$, we can use the boundedness of the probability
density of the standard normal as well as a second-order Taylor
expansion, to get for $\forall\iota>0$ and $n$ sufficiently high
\begin{eqnarray*}
\label{appsv51a}
&& \Phi \biggl(\frac{\tau+\tau\tilde{\xi}_{i,j}^n(1)-\tau\tilde{\xi
}_{i,j}^n(2)}{\xi_{i,j}^n(4)} \biggr)
\\[-1pt]
&&\qquad{} - \Phi \bigl(\tau+\tau
\tilde{\xi}_{i,j}^n(1)-\tau\tilde{\xi}_{i,j}^n(2) \bigr)
= b_i^n(1)+b_i^n(2)+b_i^n(3),
\\[-1pt]
&& b_i^n(1) = \biggl\{ \Phi \biggl(\frac{\tau+\tau\tilde{\xi
}_{i,j}^n(1)-\tau\tilde{\xi}_{i,j}^n(2)}{\xi_{i,j}^n(4)}
\biggr)
\\[-1pt]
&&\hspace*{42pt}{}- \Phi \bigl(\tau+\tau\tilde{\xi}_{i,j}^n(1)-
\tau\tilde{\xi }_{i,j}^n(2) \bigr) \biggr
\}1_{ \{|\xi_{i,j}^n(4)-1|\geq
(k_n/n)^{1/2-\iota} \}},
\\[-1pt]
&& b_i^n(2) = \Phi' \bigl(\tau+\tau
\tilde{\xi}_{i,j}^n(1)-\tau \tilde{\xi}_{i,j}^n(2) \bigr)
\\[-1pt]
&&\hspace*{35pt}{}\times \bigl(\tau+\tau\tilde{\xi
}_{i,j}^n(1)-\tau\tilde{\xi}_{i,j}^n(2)
\bigr) \bigl(\xi_{i,j}^n(4)-1\bigr)1_{ \{|\xi_{i,j}^n(4)-1|<  (k_n/n )^{1/2-\iota} \}},
\\[-1pt]
&& \bigl|b_i^n(3)\bigr|\leq K\frac{|\tau+\tau\tilde{\xi}_{i,j}^n(1)-\tau
\tilde{\xi}_{i,j}^n(2)|^2}{(1-(k_n/n)^{1/2-\iota})^3}\bigl|\xi
_{i,j}^n(4)-1\bigr|^2.
\end{eqnarray*}
For $b_i^n(1)$ and $b_i^n(3)$, we have
%
%
\[
\mathbb{E}\bigl(\bigl|b_i^n(1)\bigr|+\bigl|b_i^n(3)\bigr|
\bigr)\leq K\bigl(\tau^2\vee1\bigr)\frac
{k_n}{n}.
\]
For $b_i^n(2)$, by an application of the H\"{o}lder inequality, we
first have
%
%
\[
\mathbb{E}\bigl\llvert b_i^n(2) - \Phi'
(\tau )\tau\bigl(\xi _{i,j}^n(4)-1\bigr)1_{ \{|\xi_{i,j}^n(4)-1|<  (k_n/n)^{1/2-\iota} \}}
\bigr\rrvert \leq K|\tau|\frac{1}{\sqrt {n}}.
\]
Then
%
%
\begin{eqnarray*}
&& \mathbb{E} \Biggl( \frac{1}{\lfloor n/k_n\rfloor m_n}\sum_{j=1}^{\lfloor
n/k_n\rfloor}
\sum_{i=(j-1)k_n+1}^{(j-1)k_n+m_n} \Phi' (\tau )
\tau\bigl(\xi_{i,j}^n(4)-1\bigr)1_{ \{|\xi_{i,j}^n(4)-1|<  (k_n/n )^{1/2-\iota} \}}
\Biggr)^2
\\
&&\qquad \leq K\frac{k_nm_n}{n^2}.
\end{eqnarray*}
Therefore, altogether we get
%
%
\begin{equation}
\label{appsv51b} \mathbb{E}\bigl|A_5^n\bigr|\leq K\bigl(|\tau|\vee
\tau^2\bigr)\frac{k_n}{n}.
\end{equation}
We turn now to $A_1^n$ and $A_2^n$. Using secon-order Taylor expansion,
we can extract the leading terms in $A_n^2$. In particular, we denote
%
%
\[
\label{appsv52} \cases{ \displaystyle A_2^n(1) =
\frac{1}{\lfloor n/k_n\rfloor}\sum_{j=1}^{\lfloor n/k_n\rfloor}
\Phi'(\tau)\tau\bar{\xi}_{j}^n(1),
\vspace*{5pt}
\cr
\displaystyle A_2^n(2) =
\frac{1}{\lfloor n/k_n\rfloor}\sum_{j=1}^{\lfloor n/k_n\rfloor}\bigl(0.5
\Phi''(\tau)\tau^2\bigl(\bar{\xi
}_{j}^n(1)\bigr)^2-\Phi'(\tau)
\tau\bar{\xi}_{j}^n(2)\bigr).}
\]
With this notation, using the bounds in (\ref{appsv34}), as well as
the boundedness of $\Phi'''$, we have
%
%
\begin{equation}
\label{appsv53} \qquad\mathbb{E}\bigl|A_2^n-A_2^n(1)-A_2^n(2)\bigr|
\leq K\bigl(|\tau|^3\vee|\tau |^2\bigr) \biggl[ \biggl(
\frac{k_n}{n} \biggr)^{3/2}\vee \biggl(
\frac{1}{k_n} \biggr)^{3/2} \biggr].
\end{equation}
Further, upon denoting with $\widehat{A}_2^n(1)$ and $\widehat
{A}_2^n(2)$ the counterparts of $A_2^n(1)$ and $A_2^n(2)$ with
$\bar{\xi}_{j}^n(1)$ and $\bar{\xi}_{j}^n(2)$ replaced with
$\hat{\xi}_{j}^n(1)$ and $\hat{\xi}_{j}^n(2)$, respectively, we
have using the bounds in (\ref{appsv32}) [as well as the restriction
on the rate of growth of $k_n$ in (\ref{eqsv1})]
%
%
\begin{equation}
\label{appsv54} \qquad\mathbb{E}\bigl|A_2^n(1)+A_2^n(2)-
\widehat{A}_2^n(1)-\widehat {A}_2^n(2)\bigr|
\leq K \bigl(|\tau|\vee\tau^2\bigr) \biggl(\frac{1}{\sqrt{n}}\vee
\frac
{k_n}{n} \biggr).
\end{equation}
Thus we are left with the terms $A_1^n$, $\widehat{A}_2^n(1)$ and
$\widehat{A}_2^n(2)$. For $\widehat{A}_2^n(2)$, using
%
%
\begin{eqnarray*}\label{appsv55}
\mathbb{E}_{(j-1)k_n}^n \bigl(\hat{\xi}_{j}^n(1) \bigr)^2 &=& 2\mathbb
{E}_{(j-1)k_n}^n\bigl(\hat{\xi}_{j}^n(2)
\bigr)
\\
&=& \frac{1}{4}\frac
{1}{k_n} \biggl( \biggl(\frac{\pi}{2}
\biggr)^2+\pi-3 \biggr)+o \biggl(\frac
{1}{k_n} \biggr),
\end{eqnarray*}
we have
%
%
\begin{equation}
\label{appsv56} k_n\widehat{A}_2^n(2)
\stackrel{\mathbb{P}} {\longrightarrow} \frac
{\tau
^2\Phi''(\tau)-\tau\Phi'(\tau)}{8} \biggl( \biggl(
\frac{\pi
}{2} \biggr)^2+\pi-3 \biggr),
\end{equation}
locally uniformly in $\tau$. We finally will show that
%
%
\begin{equation}
\label{appsv57} \qquad\bigl(\matrix{\sqrt{\lfloor n/k_n\rfloor m_n}A_1^n
& \sqrt{\lfloor n/k_n\rfloor k_n} \widehat{A}_2^n(1)
}\bigr)\stackrel{\mathcal {L}} {\longrightarrow } \bigl(\matrix{Z_1(
\tau) & Z_2(\tau) }\bigr),
\end{equation}
locally uniformly in $\tau$. We have
%
%
\begin{eqnarray*}
\label{appsv58} \pmatrix{\displaystyle\sqrt{\lfloor n/k_n\rfloor
m_n}A_1^n \vspace*{5pt}
\cr
\displaystyle\sqrt{\lfloor
n/k_n\rfloor k_n} \widehat{A}_2^n(1)
} &=& \sum_{i=1}^{\lfloor n/k_n\rfloor k_n}
\pmatrix{\displaystyle
\zeta_i^n(1) \vspace*{5pt}
\cr
\displaystyle\frac{\Phi'(\tau)\tau}{2}\bigl(
\zeta_i^n(2)+\zeta_i^n(3)
\bigr)}+ \pmatrix{0 \vspace*{5pt}
\cr
\displaystyle\frac{\Phi'(\tau)\tau}{2}\tilde{\zeta}{}^n}
\nonumber
\end{eqnarray*}
with
%
%
\[
\label{appsv59} \zeta_i^n = \pmatrix{
\displaystyle\frac{1}{\sqrt{\lfloor n/k_n\rfloor m_n}} \bigl[1 \bigl(\sqrt{n}\Delta_i^nW
\leq\tau \bigr)-\Phi(\tau) \bigr] \vspace*{5pt}
\cr
\displaystyle\frac{1}{\sqrt{\lfloor n/k_n\rfloor k_n}}
\frac{\pi}{2}\bigl|\sqrt{n}\Delta _{i-1}^nW\bigr| \biggl(\bigl|\sqrt{n}
\Delta_i^nW\bigr|-\sqrt{\frac{2}{\pi}} \biggr)
\vspace*{5pt}
\cr
\displaystyle\frac{1}{\sqrt{\lfloor n/k_n\rfloor k_n}}\sqrt{\frac{\pi}{2}} \biggl( \bigl|\sqrt {n}
\Delta _i^nW\bigr|-\sqrt{\frac{2}{\pi}} \biggr) },\qquad
i\in I^n,
\]
where $I^n = \{i=(j-1)k_n+1,\ldots,(j-1)k_n+m_n, j=1,\ldots,\lfloor
n/k_n\rfloor\}$, and for $i=1,\ldots,n\setminus I^n$, $\zeta_i^n$ is
exactly as above with only the first element being replaced with zero,
and finally
%
%
\begin{eqnarray*}
\tilde{\zeta}{}^n &=& -\frac{(\pi/2)}{\sqrt{\lfloor n/k_n\rfloor
k_n}}\sum
_{j=1}^{\lfloor n/k_n\rfloor} \biggl[\bigl|\sqrt{n}\Delta
_{(j-1)k_n}^nW\bigr| \biggl(\bigl|\sqrt{n}\Delta_{(j-1)k_n+1}^nW\bigr|-
\sqrt{\frac
{2}{\pi
}} \biggr)
\\
&&\hspace*{162pt}{}          + \sqrt{\frac{2}{\pi}} \biggl(\bigl|\sqrt{n}
\Delta_{jk_n}^nW\bigr| -\sqrt {\frac
{2}{\pi}} \biggr)
\biggr],
\nonumber
\end{eqnarray*}
where we set $\Delta_0^nW=0$. With this notation, we have
%
%
\[
\mathbb{E} \bigl(\tilde{\zeta}{}^n \bigr)^2\leq
\frac{K}{k_n}.
\]
Further,
\begin{eqnarray*}
\label{appsv60} \mathbb{E}_{i-1}^n\bigl(
\zeta_i^n\bigr) &=& 0,
\\
\sum_{i=1}^{\lfloor n/k_n\rfloor
k_n}\mathbb{E}_{i-1}^n\bigl\|\zeta_i^n\bigr\|^{2+\iota}&\rightarrow& 0\qquad \forall\iota>0,
\\
\label{appsv61} \sum_{i=1}^{\lfloor n/k_n\rfloor k_n}
\mathbb{E}_{i-1}^n \bigl[\zeta _i^n
\bigl(\zeta_i^n\bigr)' \bigr]
&\rightarrow& \pmatrix{\displaystyle\Phi(\tau) \bigl(1-\Phi(\tau)\bigr) & 0 & 0
\vspace*{6pt}\cr
0 &
\displaystyle\biggl(\frac{\pi}{2} \biggr)^2 \biggl(1-\frac{2}{\pi}\biggr) & \displaystyle\frac{\pi}{2} \biggl(1-\frac{2}{\pi} \biggr)
\vspace*{6pt}\cr
0 &
\displaystyle\frac{\pi}{2} \biggl(1-\frac{2}{\pi} \biggr) & \displaystyle\frac{\pi}{2}\biggl(1-\frac{2}{\pi} \biggr)},
\end{eqnarray*}
because recall $m_n/k_n\rightarrow0$.
Combining the last two results, we have the convergence in (\ref{appsv57}), pointwise in $\tau$, by an application of
Theorem~VIII.3.32 in \cite{JS}. Application of Theorem 12.3 in \cite{Billingsley}, extends the convergence to local uniform in~$\tau$.

Altogether, the limit behavior of $\widehat{G}_n(\tau) - \Phi(\tau
)$ is
completely characterized by the limits in (\ref{appsv56})--(\ref
{appsv57}) and
%
%
\begin{equation}
\label{appsv61a} \sup_{\tau\in\mathcal{A}}\bigl|\widehat{G}_n(\tau) -
\Phi(\tau)- A_1^n - \widehat{A}_2^n(1)
- \widehat{A}_2^n(2)\bigr|= o_p \biggl(
\frac
{1}{k_n} \biggr),
\end{equation}
where $\mathcal{A}$ is a compact subset of $(-\infty,0)$, with the
result in (\ref{appsv61a}) following from the bounds on the order of
magnitude derived above.

\subsubsection{The difference \texorpdfstring{$\widehat{F}_n(\tau)-\widetilde{F}_n(\tau)$}
{Fn(tau) - Fn(tau)}}
To analyze the difference $\widehat{F}_n(\tau)-\widetilde{F}_n(\tau)$,
we use the following inequality:
%
%
\begin{eqnarray*}
\label{appsv62} &&\mathbb{P} \Bigl(\bigl|\Delta_i^n X\bigr|>
\alpha\sqrt{\widehat {V}{}^n_{j}}n^{-\varpi}
\Bigr)
\\
&&\qquad\leq\mathbb{P} \biggl(\biggl\llvert \frac
{\sqrt{\widehat{V}{}^n_{j}}}{\sigma_{((j-1)k_n)/n}}-1\biggr\rrvert >0.5
\biggr)+\mathbb{P} \bigl(\bigl|\Delta_i^n X\bigr|> 0.5\alpha\sigma
_{((j-1)k_n)/n}n^{-\varpi} \bigr).
\nonumber
\end{eqnarray*}
For the first probability on the right-hand side of the above
inequality we can use the bounds in (\ref{appsv31}), (\ref
{appsv32}) and (\ref{appsv36}), while for the second one we can use
the exponential inequality for continuous martingales with bounded
variation (see, e.g., \cite{RY}), as well as the algebraic inequality
$|\sum_ia_i|^p\leq\sum_i|a_i|^p$ for $p\in(0,1]$, to conclude
%
%
\begin{equation}
\label{appsv63} \mathbb{P} \Bigl(\bigl|\Delta_i^n X\bigr|> \alpha
\sqrt{\widehat {V}{}^n_{j}}n^{-\varpi} \Bigr)
\leq K \biggl[\frac{k_n}{n}\vee n^{-1+\iota
\varpi} \biggr]\qquad \forall\iota>0.
\end{equation}
Since $k_n/\sqrt{n}\rightarrow0$ and from the result of the previous
two subsections $\widetilde{F}_n(\tau)-\Phi(\tau) = O_p (\frac
{1}{k_n} )$, we get from here
%
%
\begin{equation}
\label{appsv64} \sup_{\tau\in\mathcal{A}}\bigl|\widehat{F}_n(\tau)-
\widetilde {F}_n(\tau)\bigr| = o_p \biggl(\frac{1}{k_n}
\biggr)
\end{equation}
for any compact subset $\mathcal{A}$ of $(-\infty,0)$.

\subsection{Proof of Theorem~\texorpdfstring{\protect\ref{teotv}}{4}}
The proof follows exactly the same steps as the proof of Theorem~\ref{teosv}, and we use analogous notation as in that proof. The only
nontrivial difference in analyzing the term $\widetilde{F}'_n(\tau) -
\widehat{G}'_n(\tau)$ regards the difference $|\widehat{C}{}^n_j-\dot
{C}^n_j|$ (and $|\widehat{C}{}^n_j(i)-\dot{C}^n_j(i)|$). For this, we
make use of the following algebraic inequality:
%
%
\begin{eqnarray*}
\label{apptv1}
\bigl|x^21_{\{|x|\leq a\}}-y^21_{\{|y|\leq a\}}\bigr|
&\leq& |x-y|^21_{\{|x-y|\leq
2a\}}+2a|x-y|1_{\{|x-y|\leq2a\}}
\\
&&{} +2|y|^21_{\{|y|>a/2\}}+a^21_{\{
|x-y|>a/2\}}.
\end{eqnarray*}
Using the above inequality, the bound in (\ref{appsv15}), as well as
the exponential inequality for continuous martingales with bounded
variation (see e.g., \cite{RY}), we have
%
%
\begin{equation}
\label{apptv2} \mathbb{P} \bigl(\bigl|\widehat{C}{}^n_j-
\dot{C}^n_j\bigr|\geq\eta_n \bigr)\leq K
\frac{k_n}{n}.
\end{equation}
Then, upon picking $\eta_n\propto n^{-q-\iota}$ for $\iota\in(0,
1/2-q)$ sufficiently small, we get $\sup_{\tau\in\mathcal
{A}}|\widetilde{F}'_n(\tau) - \widehat{G}'_n(\tau)| = o_p
(\frac
{1}{k_n} )$ for any compact subset $\mathcal{A}$ of $(-\infty,0)$.

Further, for $\widehat{G}'_n(\tau) - \Phi(\tau)$ the only difference
from the analysis of the corresponding term in the proof of
Theorem~\ref{teosv} is that now we have
\[
k_n\widehat{A}_2^n(2)\stackrel{\mathbb{P}}
{\longrightarrow} \frac
{\tau
^2\Phi''(\tau)-\tau\Phi'(\tau)}{4}
\]
and further now
%
%
\[
\label{apptv3} \pmatrix{ \sqrt{\lfloor n/k_n\rfloor
m_n}A_1^n
\vspace*{5pt}\cr
\sqrt{\lfloor
n/k_n\rfloor k_n} \widehat{A}_2^n(1)}
= \sum_{i=1}^{\lfloor n/k_n\rfloor k_n}
\pmatrix{\displaystyle \zeta_i^n(1)
\vspace*{5pt}\cr
\displaystyle \frac{\Phi'(\tau)\tau}{2}\zeta_i^n(2)}
\]
with
%
%
\begin{eqnarray*}\label{apptv4}
\bigl(\zeta_i^n\bigr)' &=&
\biggl(\matrix{\displaystyle\frac{1}{\sqrt{\lfloor n/k_n\rfloor m_n}} \bigl[1 \bigl(\sqrt{n}\Delta_i^nW
\leq\tau \bigr)-\Phi(\tau) \bigr]}
\\
&&\hspace*{48pt} \matrix{\displaystyle\frac{1}{\sqrt {\lfloor n/k_n\rfloor k_n}}\bigl(\bigl(\sqrt{n}
\Delta_i^nW\bigr)^2-1\bigr)}\biggr),\qquad
i\in I^n,
\end{eqnarray*}
where $I^n = \{i=(j-1)k_n+1,\ldots,(j-1)k_n+m_n, j=1,\ldots,\lfloor
n/k_n\rfloor\}$, and for $i=1,\ldots,n\setminus I^n$, $\zeta_i^n$ is
exactly as above with only the first element being replaced with zero.
From here the analysis of $\widehat{G}'_n(\tau) - \Phi(\tau)$ is done
exactly as that of the corresponding term in the proof of Theorem~\ref{teosv}.

We are left with showing the result in the case when jumps in $X$ can
be of infinite activity (under the conditions in the theorem). We again
follow the steps of the proof of Theorem~\ref{teosv}. We replace
$R_j^{(4)}$ with $\widehat{C}{}^n_j-\dot{C}^n_j$ in $\bar{\xi
}_j^n(1)$ and $\bar{\xi}_j^n(2)$ and similarly we replace
$\widetilde{R}_{i,j}^{(4)}(i)$ with $\widehat{C}{}^n_j-\dot{C}^n_j -
(\Delta_i^nX)^21_{\{|\Delta_i^nX|\leq\alpha n^{-\varpi}\}}+|\Delta
_i^nA+\Delta_i^nB|^2$ in $\tilde{\xi}_{i,j}^n(1)$ and
$\tilde{\xi}_{i,j}^n(2)$.

Using the inequality in (\ref{apptv1}), and since $\int_{E}|\delta
^{Y}(x)|^{\beta'}\nu(dx)<\infty$ (upon localization that bounds the
size of the jumps), we have
%
%
\begin{eqnarray}
\label{apptv5} \qquad\mathbb{E}\bigl\llvert \bigl(\Delta_i^nX
\bigr)^21_{\{|\Delta_i^nX|\leq\alpha
n^{-\varpi }\}}-\bigl(\Delta_i^nA+
\Delta_i^nB\bigr)^2\bigr\rrvert
^p\leq Kn^{-1-(2p-\beta
')\varpi}
\nonumber\\[-10pt]\\[-10pt]
\eqntext{\forall p\geq\beta'/2}
\end{eqnarray}
and from here
%
%
\begin{equation}
\label{apptv6} \mathbb{E}\bigl|\widehat{C}{}^n_j-
\dot{C}^n_j\bigr|^p\leq Kn^{p-1-(2p-\beta
')\varpi
}\qquad \forall p\geq1.
\end{equation}
Using the bounds in (\ref{apptv5}) and (\ref{apptv6}), we can prove
exactly as in the proof of Theorem~\ref{teosv} for some deterministic
sequence of positive numbers $\eta_n$
%
%
\begin{eqnarray*}
\label{apptv7} &&\mathbb{P} \bigl(\bigl(\bigl|\chi_{i,j}^n(1)\bigr|+\bigl|
\chi_{i,j}^n(2)\bigr|\bigr)>\eta_n \bigr)
\\
&&\qquad\leq K \biggl[\frac{1}{\eta_n^{p}[n^{p/2}\wedge(n/k_n)^p\wedge
k_n^{3p/2}]}\vee \biggl(\frac{k_n}{n}
\biggr)^{1/(1+\iota)}\frac{1}{\eta_n^{\iota}}
\\
&&\hspace*{123pt}{}\vee\frac{n^{-1+\beta'/2}}{\eta_n^{\beta
'}}\vee\frac{n^{-(2-\beta')\varpi}}{\eta_n^{1/2}\wedge\eta
_n\sqrt {k_n}} \biggr].
\end{eqnarray*}
From here, using the rate of growth condition in (\ref{eqtv4}), upon
appropriately choosing~$\eta_n$, we get
%
%
\begin{equation}
\label{apptv8} \sup_{\tau\in\mathcal{A}}\bigl|\widetilde{F}'_n(
\tau) - \widehat {G}'_n(\tau )\bigr| = o_p
\biggl(\frac{1}{k_n} \biggr)
\end{equation}
for any compact subset $\mathcal{A}$ of $(-\infty,0)$.

We turn next to $\widehat{G}'_n(\tau) - \Phi(\tau)$, and we derive the
bounds of those terms in the decomposition of the latter which are
different from the case of finite jump activity proved above (the term
$A_5$ is identically zero since $\sigma_t$ is constant). First, for
$A_3^n$, using (\ref{apptv6}) as well as the independence of $W_t$
and $Y_t$, we have
%
%
\begin{eqnarray}
\label{apptv9} \qquad\qquad A_3^n &\leq& K \bigl(\sqrt{|\tau|}\vee
\tau^2\bigr)
\nonumber\\[-8pt]\\[-8pt]
&&{}\times \biggl(\frac{1}{\sqrt {\lfloor n/k_n\rfloor m_n}}n^{1/2-(4-\beta')\varpi/2}\vee
k_n^{1/4}n^{-1/4-(4-\beta')\varpi/4}\vee
\frac{1}{\sqrt {n}} \biggr).\nonumber
\end{eqnarray}
Next, if we exclude $\widehat{C}{}^n_j-\dot{C}^n_j$ from $\bar {\xi
}_j^n(1)$ and $\bar{\xi}_j^n(2)$, we get for $A_4$, using (\ref
{apptv5}) and (\ref{apptv6}), as well as applying the H\"{o}lder
inequality,
%
%
\begin{eqnarray}
\label{apptv10}
\mathbb{E}\bigl|A_4^n\bigr|
&\leq& K\bigl(|\tau|\vee
\tau^2\bigr)
\nonumber\\[-8pt]\\[-8pt]
&&\times{} \biggl(\frac{1}{k_n^{3/2}}\vee
\frac{1}{\sqrt{\lfloor n/k_n\rfloor
}k_n}\vee \frac{n^{-(2-\beta')\varpi}}{\sqrt{k_n}}\vee
n^{1-(4-\beta
')\varpi} \biggr).\nonumber
\end{eqnarray}
Combining the bounds in (\ref{apptv9})--(\ref{apptv10}), and taking
into account the growth condition in (\ref{eqtv4}), we get
%
%
\begin{equation}
\label{apptv12} \sup_{\tau\in\mathcal{A}}\bigl|\widehat{G}'_n(
\tau) - \Phi(\tau)- A_1^n - \widehat{A}_2^n(1)
- \widehat{A}_2^n(2)\bigr|= o_p \biggl(
\frac
{1}{k_n} \biggr),
\end{equation}
where $\mathcal{A}$ is a compact subset of $(-\infty,0)$. The limit
behavior of the triple $(A_1^n,\widehat{A}_2^n(1),\widehat{A}_2^n(2))$ is derived as in the finite jump activity case in the
first part of the proof and this together with (\ref{apptv8}) and
(\ref{apptv12}) yields the stated result in the case of infinite
variation jumps.

\section*{Acknowledgments}
We would like to
thank Dobrislav Dobrev, Jean Jacod, Per Mykland, Mark Podolskij, Markus
Reiss, Mathieu Rosenbaum and many seminar participants for helpful
comments and suggestions. We also thank an Associate Editor and a
referee for careful read and many constructive comments.




\printaddresses

\end{document}